\documentclass{article}

\usepackage[numbers,sort&compress]{natbib}

\usepackage{mathtools}
\usepackage{authblk}
\usepackage[title]{appendix}
\usepackage{geometry}
\usepackage{subfiles}
\usepackage{caption}
\usepackage{subcaption}
\usepackage{setspace}
\usepackage{amsthm}
\usepackage{amsmath}
\usepackage{amssymb}
\usepackage{float}
\usepackage{fancyhdr}
\usepackage{graphicx}
\usepackage[colorlinks,            linkcolor=black,            anchorcolor=black,            citecolor=black            ]{hyperref}
\usepackage{listings}
\usepackage{tikz}
\usetikzlibrary{arrows,decorations.pathmorphing,backgrounds,positioning,fit,petri,arrows.meta,bending,positioning}

\newtheorem{corollary}{Corollary}

\newtheorem{con-proposition}{``Proposition''}
\newtheorem{proposition}{Proposition}
\newtheorem{theorem}{Theorem}
\newtheorem{lemma}{Lemma}

\newtheorem{assumption}{Assumption}

\newcommand{\Div}{\nabla\cdot}
\newcommand{\p}{\partial}

\providecommand{\keywords}[1]
{
	\small	
	\textbf{{Keywords:}} #1
}
\providecommand{\msc}[1]
{
	\small	
	\textbf{{Mathematics Subject Classification:}} #1
}

\newcommand{\hl}[1]{{#1}}

\title{A tumor growth model with autophagy: the reaction-(cross-)diffusion system and its free boundary limit }

\author{Xu'an Dou\thanks{School of Mathematics Science, Peking University, Beijing, 100871, China (dxa@pku.edu.cn)} \hspace{2cm}{Jian-Guo Liu}\thanks{Department of Mathematics and Department of Physics, Duke University, Box 90320, Durham, NC27708, USA(jliu@phy.duke.edu).}\hspace{2cm}
{Zhennan Zhou }\thanks{Beijing International Center for Mathematical Research, Peking University, Beijing, 100871, China (zhennan@bicmr.pku.edu.cn).}}

\begin{document}
\maketitle
\begin{abstract}
In this paper, we propose a tumor growth model to incorporate and investigate the spatial effects of autophagy. The cells are classified into two phases: normal cells and autophagic cells, whose dynamics are also coupled with the nutrients. First, we construct a reaction-(cross-)diffusion system describing the evolution of cell densities, where the drift is determined by the negative gradient of the joint pressure, and the reaction terms manifest the unique mechanism of autophagy. Next, in the incompressible limit, such a cell density model naturally connects to a free boundary system, describing the geometric motion of the tumor region. Analyzing the free boundary model in a special case, we show that the ratio of the two phases of cells exponentially converges to a ``well-mixed" limit. Within this ``well-mixed" limit, we obtain an analytical solution of the free boundary system which indicates the exponential growth of the tumor size in the presence of autophagy in contrast to the linear growth without it. Numerical simulations are also provided to illustrate the analytical properties and to explore more scenarios. 

\end{abstract}


\msc{35Q92, 35R35, 76D27, 92-10, 92D25}


\section{Introduction}\hl{
 Autophagy, coined from ``self-eating'' in Greek, is a (catabolic) process in cells, during which cells degrade its own constituents. The first-discovered and most fundamental function of autophagy is to help cells survive when the nutrient is lacking. This response to nutrient depletion plays an important role in tumor progression (e.g. reviews \cite{crlevine2007autophagy,crwhite2012deconvoluting,crlorin2013autophagy,crli2020autophagy}). On the other hand, tumor growth is a complicated  dynamical process which describes the mechanical motion of collective tumor cells, while the subtle evolution of each cell also plays a significant role. In this work, we consider a macroscopic cell density model to investigate the autophagy effect on tumor spatial growth. The reaction part of the model is related to the microscopic behavior of individual cells, which are based on the following three scientific features of autophagy.}



\hl{
{First, autophagy induction and termination are regulated by nutrients.} Autophagy is drastically induced when there is shortage of nutrients (amino acids, glucose, oxygen ...), while it is inhibited in nutrient-rich condition \cite{ohsumi2014historical,russell2014autophagy}.  For example, a series of experiments have shown that, lack of amino acids induces autophagy to regenerate nutrients, while the restoration of amino acids terminates autophagy \cite{yu2010autophagy}.}
\hl{
{Second, autophagy provides nutrients.} When autophagy is induced, cells degrade its own cytoplasm and organelles, breaking down the macromolecules (such as proteins) into monomer units (such as amid acids), therefore generating nutrients for reuse \cite{mizushima2007autophagy_process_and_function}. This supply is crucial for cells to survive under nutrient deprivation, and many experiments have shown that autophagy-deficient cells rapidly lose their viability under nutrient starvation \cite{tsukada1993isolation,komatsu2005impairment,mizushima2011autophagy}.} 
\hl{
{Third, autophagy inhibits growth and promotes death.} Autophagy is often correlated with the suppression of cell growth, since autophagy is a degradation process and cell growth needs synthesis of macromolecules. Moreover, it is revealed by various experiments that autophagy has a causal inhibition effect on cell growth \cite{neufeld2012autophagy}. Besides, a cell under autophagy may ``eat itself to die'' due to over-consumption of its own constituents (cellular organelles and cytoplasmic content) \cite{bialik2018autophagy}. Combining these effects, at the macroscopic level we assume autophagy has a negative effect on the net growth rate of cells.}

\hl{Based on the above three characteristics, we construct the following reaction system:}
\begin{spacing}{1.3}\begin{equation}\label{model:ODE}
	\begin{cases}
	\frac{dn_1}{dt}=G(c)n_1-K_1(c)n_1+K_2(c)n_2.\\
	\frac{dn_2}{dt}=(G(c)-D)n_2+K_1(c)n_1-K_2(c)n_2.\\
	\frac{dc}{dt}=-\lambda(t)(c-c_B(t))-\psi(c)(n_1+n_2)+an_2.\\
	\end{cases}
	\end{equation} 
\end{spacing}

\tikzset{global scale/.style={
		scale=#1,
		every node/.append style={scale=#1}
}}
\begin{figure}[H]
	\centering	\begin{tikzpicture}
	[global scale=0.7,auto,inner sep=2mm,density/.style={draw=black!50,fill=black!20,minimum height=20mm,minimum width=60mm},nutrient/.style={rectangle,draw=green!50,fill=green!20,thick,minimum size=30mm}]
	\node[density,font=\fontsize{22}{18}\selectfont](n1) at (2,6) {Normal Cell: $n_1$}
	edge [->,looseness=4,in=170,out=190,-{Latex[]}] node[font=\fontsize{15}{15}\selectfont] {$G$} (n1);
	\node[density,font=\fontsize{22}{15}\selectfont](n2) at (2,2){Autophagic Cell: $n_2$}
	edge [->,looseness=4,in=170,out=190,-{Latex[]}] node[font=\fontsize{15}{15}\selectfont] {$G$} (n2);
	\node[nutrient,font=\fontsize{22}{15}\selectfont](c) at (10,4){Nutrient: $c$};
	\draw [->,dashed,-{Latex[]}](c) to node[font=\fontsize{15}{15}\selectfont]{$\psi$}(n1);
	\draw [->,dashed,-{Latex[]}](c) to node[font=\fontsize{15}{15}\selectfont]{$\psi$}(n2);
	\draw [->,dashed,-{Latex[]}](n2) to[in=240,out= 0] node[swap,font=\fontsize{15}{15}\selectfont]{\textcolor{red}{$a$}}(c);
	\draw [->,-{Latex[]}](n1) to[in=75,out=285] node[font=\fontsize{15}{15}\selectfont]{$K_1$}(n2);
	\draw [->,-{Latex[]}](n2) to[out=105,in=255] node[font=\fontsize{15}{15}\selectfont]{$K_2$}(n1);
	\node(outside) [above right=of c]{}
	edge [->,-{Latex[]}] node[font=\fontsize{15}{15}\selectfont] {$\lambda c_B$}(c);
	\node(outside21) at (10,1.1){};
	\node(outside22) [below right=of outside21]{}
	edge [<-,{Latex[]}-] node[font=\fontsize{15}{15}\selectfont] {$\lambda c$}(c);
	\node(dead n2) [below left= of n2]{}
	edge [<-,{Latex[]}-] node[swap,font=\fontsize{15}{15}\selectfont] {\textcolor{red}{$D$}}(n2);
	\end{tikzpicture}
	\caption{Illustration of the ODE model. The two gray squares denote the densities of normal cells $n_1$ and autophagic cells $n_2$. And the green square denotes the nutrient concentration $c$. Normal cells change into autophagic cells with a transition rate $K_1$. And autophagic cells change into normal cells with a transition rate $K_2$. Normal cells and autophagic cells both grow with a net growth rate $G$. Dashed lines show the consumption or supply of nutrients by cells. Both normal cells and autophagic cells consume nutrients with a rate $\psi$. The key assumption is that autophagic cells will ``kill'' themselves with an extra death rate $D$ to provide nutrients with a supply rate $a$. Nutrients are added with a rate $\lambda c_B$ and discharged with a rate $\lambda c$. }
\end{figure}
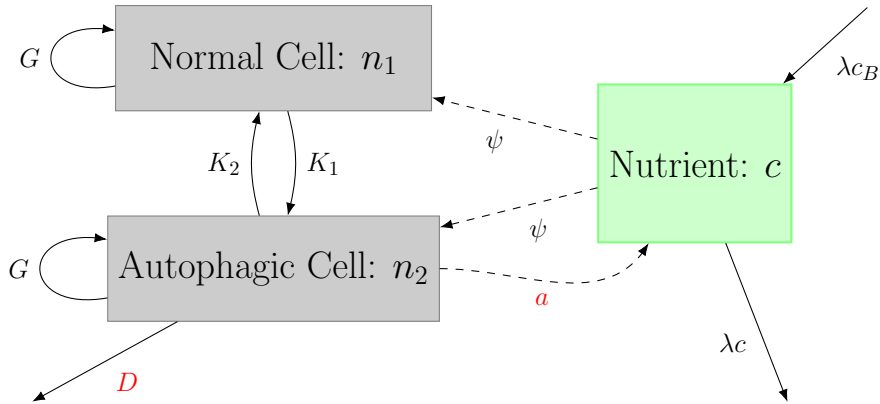

In this ODE model \eqref{model:ODE}, cells are classified into two phases: normal cells and autophagic cells. Their densities are $n_1,n_2$ while $c$ denotes the concentration of nutrients. 
$G(c)$ is the net growth rate of normal cells, which could be understood as the birth rate minus the death rate, and is increasing with respect to the concentration of nutrients. Autophagy is induced under nutrient depletion, and is terminated due to the restoration of nutrients. Therefore, we assume that normal cells are transformed into the autophagic phase with a rate $K_1(c)$, which is decreasingly with respect to $c$, and that autophagic cells turn into normal phase with a rate $K_2(c)$, which is increasing with respect to $c$. That is, we have the following assumptions on $G,K_1,K_2$:
\begin{equation}\label{ap-2GK}
G'(c)\geq0;\quad K_1(c),\,K_2(c)>0,\,\,\text{and} \,\, K_1'(c)\leq0,\,K_2'(c)\geq0.
\end{equation}

One of the key ingredients of this model is that the tumor cells are categorized into two phases connecting by a reversible process. \hl{Such a representation is in fact a manifestation of the regulation mechanism of the autophagy-related (Atg) proteins at the microscopic level.  In the absence of adequate nutrient,  conductor proteins change the mode of the pre-autophagosomal structure (PAS, an assemble of proteins related in autophagy) from the Cvt pathway (which appears in a nutrient-rich condition) to autophagy (experiments \cite{kawamata2008organization,cheong2008atg1} and the review \cite{nakatogawa2009dynamicsteast}), and the change in the PAS indicates that a cell changes from the normal phase to the autophagic phase.}

\hl{At the macroscopic level, we have also taken the assumptions that autophagic cells provide nutrients proportional to its density with a rate $a$, and that autophagy has a negative effect $D$ on growth rate. These are direct implications of the second and the third features of autophagy, respectively.} In this work we assume
\begin{equation}\label{ap-1Da}
D>0,\quad a>0,
\end{equation} are constants for simplicity.

The rate of change of nutrient concentration $c$ is determined by the influence of the background concentration via the term $\lambda(c_B-c)$, the consumption by both types of cells and a supply from the degradation of autophagic cells.	
For the consumption rate $\psi(c)$ we assume it is non-negative and increasing with $c$. When there is no nutrient the consumption rate shall be zero. Moreover, we assume there exists a critical nutrient concentration $c_0$, such that $\psi(c_0)=a$. To summarize, we have the following assumption on the  consumption rate $\psi(c)$:
\begin{equation}\label{ap-3psi}
\psi(0)=0,\quad\psi'(c)>0,\,\, \text{and},\, \exists \, c_0>0 \, \text{s.t.}\ \psi(c_0)=a.
\end{equation}If $c>c_0$, then $\psi(c)>a$, which means the nutrient is sufficient enough that autophagic cells consume more nutrients than they supply. While when $0<c<c_0$, $\psi(c)-a<0$, which means autophagic cells supply more nutrients than they consume. 


\hl{Autophagy has been recognized as a fundamental and ubiquitous process by biologists. It is well-conserved in yeasts, mammals and plants. Besides its most typical mode -- response to starvation, it can also be induced by various stimuli such as damaged organelles,  pathogens, high temperature and hormones and it is related with a wide range of processes such as cancer, heart disease, aging and so on  \cite{crlevine2007autophagy,klionsky2007autophagy,klionsky2020autophagy}. Its role is often complicated. For example in cancer, while autophagy can help tumor cells survival under nutrient deprivation, via degrading damaged proteins and organelles, autophagy can prevent tumor initiation in the early stage \cite{crlevine2007autophagy,crwhite2012deconvoluting,crlorin2013autophagy,crli2020autophagy}.  }

\hl{
In this work, we focus on a major function of autophagy -- response to starvation, in a particular scenario -- tumor growth. As a clump of tumor grows, its central part becomes far from the vascular supply, therefore facing a lack of nutrients  \cite{roose2007mathematical}. Before the tumor cells develop its own vascular (in the avascular stage), the nutrient availability varies strongly in different parts of the tumor -- it is sufficient near the boundary but may be insufficient in the center. Autophagy is crucial for cells in the center to survive nutrient starvation \cite{degenhardt2006autophagy,crlevine2007autophagy,kang2008or}.} 

Despite the importance and popularity of autophagy, at population level few mathematical models are proposed. As far as we know, \cite{jin2014mathematical-autophagy} is the first to establish an ODE model for cell population dynamics with the effects of autophagy and actually the ODE system \eqref{model:ODE} can be seen as a modified version of their model. However, as aforementioned the spatial heterogeneity is crucial in tumor growth. To our best knowledge, despite autophagy plays a crucial role in tumor growth, there is no spatial population model developed to understand its effect. The motivation of this paper is to study how autophagy influences the spatial growth of tumor. Hence, we propose and study a spatial model of tumor growth in the presence of autophagy. 

For tumor growth in space, there have been numerous continuum models in literature \cite{roose2007mathematical,lowengrub2010nonlinear,perthame2015some}, most of which fall into two categories: either they describe the densities of tumor from a fluid mechanics point of view \cite{bresch2010computational,ranft2010fluidization}, or they treat the tumor as an expanding domain $\Omega(t)$ and describe its geometric motion, where free boundary problems arise \cite{friedman2004hierarchy,friedman2007mathematical}. The connection between these two kinds of models has been established through the ``incompressible limit'' \cite{perthame2014heleAsym,ahelemellet2017hele}. 

\hl{
In the avascular stage, nutrient is often modeled as diffusing from the boundary of tumor \cite{roose2007mathematical}. Thus the nutrient is sufficient at the boundary but becomes insufficient in the center as tumor grows. This limitation of nutrient is a crucial factor of tumor growth and is responsible for many biological/math behavior such as dormant steady state, propagation speed of the boundary, necrosis \cite{https://doi.org/10.1002/sapm1972514317,PAN2018362,liu2019analysis,david:hal-02515263,guillen2020heleshaw,perthame2014traveling}, etc. Therefore it is worth exploring the interplay between autophagy and nutrient in the spatial structure of solid tumors. }

 \hl{Our spatial tumor growth models are constructed as follows.} First, we extend the ODE system \eqref{model:ODE} to a PDE system describing the evolution of cell densities, as a tumor growth model with the effects of autophagy. We denote the local density at position $x\in\mathbb{R}^d$, time $t\geq0$ by $n_1(x,t)$ (or $n_2(x,t),c(x,t)$) instead of a homogeneous density $n_1(t)$ (or $n_2(t),c(t)$). The spatial motions of cells and nutrients are different. While for nutrients we assume they diffuse quickly, the spatial motion of cells is more complicated. We assume both kinds of cells are driven by a velocity field $v$. The velocity field $v$ is induced by the local pressure $p(x,t)$ through Darcy's law: $v=-\nabla p$ from a fluid mechanical viewpoint \cite{perthame2014heleAsym}. We connects the pressure to densities via a simple constitute law $p(x,t)=\frac{\gamma}{\gamma-1}(n_1(x,t)+n_2(x,t))^{\gamma-1}(\gamma>1)$. Biologically the interpretation is that proliferating cells exert a pressure on their nearby cells, which pushes them to move. The full presentation of this PDE model is given in section \ref{subsc:PDEmodel}.

Mathematically our assumptions on the spatial motion of cells lead to a porous-media type cross-diffusion system with reactions, of which a systematic treatment is still lacking in the current literature. In the absence of nutrients, existence is first established in \cite{bertsch2010free,bertsch2012nonlinear} with a non-vacuum assumption on initial data. Recently it is shown that well-posedness can be established without this assumption \cite{carrillo2018splitting,gwiazda2019two-hyperbolic,price2020global}. Moreover, the incompressible limit to a Hele-Shaw problem has been established \cite{bubba2019-helelimit2} and the case when two kinds of cells have different
mobilities has been studied \cite{lorenzi2017interfaces-differentmobility,kim2020interface}.

The cross-diffusion system with reactions, which is also coupled with the dynamics of nutrients, is formidable for direct analysis. Such a cross-diffusion model can be interpreted as a cell density model, whose incompressible limit has been identified in previous work \cite{perthame2014heleAsym,bubba2019-helelimit2,david:hal-02515263}. Heuristically, as $\gamma\rightarrow+\infty$, the cell density model becomes ``incompressible'' and converges to a free boundary model. The resulting free boundary problem is more tractable to obtain analytical results and is a useful approximation for the cell density model with large $\gamma$. There have been fruitful results studying the tumor growth in one-species cases with the help of the incompressible limit, for example traveling solutions \cite{perthame2014traveling}, analytical solutions \cite{liu2019analysis,liu2019towards} and numerical schemes \cite{liu2018accurate,guillengonzalez2019from}. The incompressible limit in tumor growth model is first rigorously established in \cite{perthame2014heleAsym} and quite recently two-species cases \cite{bubba2019-helelimit2,degond2020incompressible} and the case with nutrients \cite{david:hal-02515263} are studied, which are more relevant to our PDE model involving autophagy.

Next, we derive a free boundary model with a formal incompressible limit from the cell density model. Here we do not justify the limit rigorously but study the limited system directly. In the free boundary model, the cells are assumed incompressible and the tumor is described by an expanding domain $\Omega(t)$. The total density $n=n_1+n_2$ inside $\Omega(t)$ is assumed to be a constant. And the pressure $p$ and the concentration of nutrients $c$ evolve as $\Omega(t)$ evolves. To describe the interplay of two kinds of cells, a density fraction $\mu(x,t)$ is introduced. $\mu(x,t)$ stands for the ratio of normal cells in all cells, corresponding to $\frac{n_1(x,t)}{n_1(x,t)+n_2(x,t)}$ in the cell density model. The formal derivation and full presentation of this system are given in subsection \ref{subsc:fbmodel}. \hl{We remark that this limit has been justified in  a subsequent work \cite{liu2021existence}.}

Free boundary models describing tumor consisting of different kinds of cells have been proposed directly, rather than from an incompressible limit, in literature. Ward and King \cite{ward1997mathematical,ward1999mathematical} consider models for living and dead cells. Pettet et al. \cite{pettet2001the} constructed a model for proliferating cells and quiescent cells. And a general free boundary framework was given by Friedman \cite{friedman2004hierarchy}. These multi-species free boundary models have attracted much mathematical interest for well-posedness \cite{cui2003a,general2003,cui2005global}, stationary solutions and stability \cite{chen2005hyperbolic,cui2008asymptotic}. We note that in the presence of nutrients, the free boundary models are not exactly the incompressible limit of the cell density models when there is necrotic core \cite{perthame2015some,perthame2014traveling}.

In summary, to study spatial effects of autophagy in tumor growth, we construct a compressible reaction-(cross-)diffusion system describing the cell density, whose incompressible limit is a free boundary model. To investigate these two systems, first we analyze the free boundary model. An interesting result here is that in a special case, we can obtain a clear asymptotic result for the density fraction $\mu$: $\mu$ converges to a spatial homogeneous constant state. We call this phenomenon the ``well-mixed'' limit. We obtain two results characterizing this limit: a uniform convergence result (Theorem \ref{thm:uni-con}) and convergence in the ``$L^{2n}(n\in\mathbb{N}^+)$ norm'' under a family of conditions (Theorem \ref{thm:L2n}).

This ``well-mixed'' limit of density fraction $\mu$ allows us to further simplify the model and derive analytical solutions with a similar manner as in \cite{liu2019analysis}. From the analytical solution we observe that when the nutrient provided by autophagy is sufficient, the tumor will expand exponentially in contrast to linear growth rate without autophagy in \cite{liu2019analysis}. The essential mechanism behind the exponential growth is that when the extra nutrient supply rate $a$ from autophagy has a bigger effect than the extra death rate $D$ due to autophagy, the tumor can access to sufficient nutrients whatever big it is. \hl{These analytical observations are in accordance with the biological observation that autophagy helps cells where the nutrient is insufficient, and therefore promotes tumor growth \cite{crwhite2012deconvoluting,kang2008or,degenhardt2006autophagy}.}

\paragraph{Arrangement of this paper}The rest of this paper is arranged as follows. We present the two systems, compressible and incompressible, and discuss the assumptions on the free boundary model in Section \ref{sc:model2}. Under these assumptions, in Section \ref{sc:math-proof} we analyze the free boundary model and prove the ``well-mixed'' limit.
Within the ``well-mixed'' limit, in Section \ref{sc:auto} we obtain an analytical solution in a special case to further investigate the model, in particular the expansion speed of the tumor. We also numerically demonstrate analytical properties and explore more scenarios.


\section{Modeling and interpretations}\label{sc:model2}
In this section, we present the two systems for tumor growth with autophagy: the compressible cell density model and the incompressible free boundary boundary model. Moreover, the model assumptions are presented with detailed interpretations. 
\subsection{Compressible cell density model}\label{subsc:PDEmodel}
We extend the ODE system \eqref{model:ODE} to involve spatial effects. We assume cells are driven by the negative gradient of pressure while nutrients diffuse quickly, which leads to the following PDE system:
\begin{spacing}{1.3}
\begin{equation}\label{model:PDE-systems}
\begin{cases}
\frac{\partial n_1}{\partial t}-\Div(n_1 \nabla p)=G(c)n_1-K_1(c)n_1+K_2(c)n_2,\quad x\in\mathbb{R}^d, t>0.\\
\frac{\partial n_2}{\partial t}-\Div(n_2 \nabla p)=(G(c)-D)n_2+K_1(c)n_1-K_2(c)n_2,\quad x\in\mathbb{R}^d, t>0.\\
\epsilon\frac{\p c}{\p t}-\Delta c+ \psi(c)(n_1+n_2)=an_2,\quad x\in\mathbb{R}^d, t>0.\\
c(x,t)\rightarrow c_B,\ |x|\rightarrow\infty.
\end{cases}
\end{equation}\end{spacing}
Here $n_1(x,t),n_2(x,t)$ are local densities of two kinds of cells and $c(x,t)$ represents the local concentration of nutrients, at position $x\in\mathbb{R}^d$, time $t>0$. As in the ODE model, $G(c)$ is the net growth rate of normal cells and $K_1(c),K_2(c)$ are transition rates with which two phases of cells change into each other. $\psi(c)$ is the consumption rate of nutrients. The key ingredient for autophagy in the ODE model is ``inherited'': $D$ is the extra death rate of autophagic cells and $a$ is the nutrient supply rate from autophagy. The parameters satisfy \eqref{ap-1Da}$\sim$\eqref{ap-3psi} as in the ODE model. 

For the spatial motion of cells, we take a fluid mechanical point of view. We assume they are driven by a velocity field which equals to the negative gradient of the pressure $p$ (Darcy's law) \cite{perthame2014heleAsym}. And the pressure arises from the mechanical contact between cells. Here we assume $p$ is a power of total density $n=n_1+n_2$, precisely:
\begin{equation}\label{eq:intro-p}
p(x,t)=\frac{\gamma}{\gamma-1}(n(x,t))^{\gamma-1},\quad n(x,t):= n_1(x,t)+n_2(x,t),\ \gamma>1.
\end{equation}

The choice \eqref{eq:intro-p} and Darcy's law lead to porous-media type diffusion. Indeed if we add the first two equations in \eqref{model:PDE-systems} and introduce the density fraction $\mu=\frac{n_1}{n_1+n_2}$, then we obtain
\begin{equation}\label{eq:porous}
\frac{\p n}{\p t}=\Delta(n^\gamma)+\mu G(c) n+(1-\mu)(G(c)-D)n,
\end{equation}where $\Delta(n^\gamma)$ term arises.

For the spatial behavior of nutrients $c$, we assume they diffuse quickly.  $\epsilon>0$ is a parameter reflecting the time scale of the evolution of nutrients, which is usually smaller than cells. In the following we consider both the quasi-static case, i.e. $\epsilon=0$, and the case $\epsilon=1$. The boundary value $c_B$ corresponds to the nutrients flux supplied by the environment. This system for nutrients has an important variant. When initially total density $n=n_1+n_2$ has a compact support, then at least formally from the property of porous-media type diffusion \eqref{eq:porous},  we deduce that $n(\cdot,t)$ has a compact support for every $t>0$ (e.g., see section 4.2 in \cite{perthame2015some}). Thus $\Omega(t):=\{x: n(x,t)>0\}$, the region occupied by tumor, is a bounded domain for all $t>0$. Therefore we could use the following model for nutrients:
\begin{equation}\label{eq:c-compact}
\begin{cases}
\epsilon\frac{\p c}{\p t}-\Delta c+\psi(c)(n_1+n_2)=an_2,\quad x\in\Omega(t),\\
c=c_B,\quad x\notin\Omega(t).\\
\Omega(t):=\{x: n(x,t)>0\}.
\end{cases}
\end{equation}The above governing equation for nutrient \eqref{eq:c-compact} is used in the free boundary model.

To complete the system \eqref{model:PDE-systems},\eqref{eq:intro-p}, we should involve initial values of $n_1,n_2$:
\begin{equation}
\begin{cases}
n_1(x,0)=n_1^0(x)\geq 0 ,\\
n_2(x,0)=n_2^0(x)\geq 0,
\end{cases}
\end{equation} and an initial value for nutrients $c$ if $\epsilon>0$.

Moreover, if we try to derive an equation for the density fraction $\mu=\frac{n_1}{n_1+n_2}$, with straightforward calculation we have (at least formally, since the definition of $\mu$ out the support of $n$ is ambiguous)
\begin{equation}\label{eq:model-mu-Com}
\frac{\partial \mu}{\partial t}-(\nabla \mu)\cdot(\nabla p)=-\mu K_1(c)+(1-\mu)K_2(c)+D\mu(1-\mu),\ 
\end{equation}which indicates the hyperbolic nature of the dynamics of two-species. This equation is used to characterize the two species dynamics instead of the equations for each species $n_1,n_2$ in the free boundary model.

\hl{We remark that the well-posedness of system \eqref{model:PDE-systems} has been established in \cite{liu2021existence}.}
	\subsection{Incompressible free boundary model and assumptions}\label{subsc:fbmodel}
	In this subsection, we derive and present the free boundary model. Then we give assumptions for the analysis in the next section.
\subsubsection{Derivation of the free boundary model: a formal incompressible limit}
	We derive the free boundary model via a formal incompressible limit from the compressible cell density model \eqref{model:PDE-systems}. The formal derivation is a direct adaptation of the one-species case \cite{perthame2014heleAsym,perthame2015some}. \hl{However the rigorous justification, which has been done in a subsequent work \cite{liu2021existence}, needs new key estimates.}

	To be concise we introduce two shorthand notations
	\begin{equation}\label{eq:g1g2}
	G_1(c):=G(c),\quad G_2(c):=G(c)-D,
	\end{equation}which are the net growth rates of normal cells and autophagic cells, respectively.
	
	 First we multiply \eqref{eq:porous} with $\gamma n^{\gamma-2}$, which is the derivative of $p$ with respect to the total density $n$, and then we obtain
	\begin{equation}\label{eq:p}
	\frac{\partial p}{\partial t}-|\nabla p|^2=(\gamma-1) p( \Delta p+\mu G_1(c)+(1-\mu)G_2(c)).
	\end{equation}
	Formally let $\gamma\rightarrow\infty$ in \eqref{eq:p}, one obtains
	\begin{equation}\label{dev-con-relation}
p( \Delta p+\mu G_1(c)+(1-\mu)G_2(c))=0.
	\end{equation}
	
	Assuming $p$ is bounded, as $\gamma\rightarrow\infty$ the relation of $p=\frac{\gamma}{\gamma-1}n^\gamma$  leads to
	\begin{equation}\label{dev-pn}
		n(x,t)\begin{cases}
		=1,\quad p(x,t)>0.\\
		\in[0,1],\quad p(x,t)=0.
		\end{cases}
	\end{equation}
	
	To avoid the formation of necrotic core (c.f \cite{perthame2014traveling,perthame2015some,cui2001analysis}), we assume the following weighted net growth rate is always non-negative:
	\begin{equation}\label{dev:as2}
	\mu G_1+(1-\mu)G_2\geq 0,\quad \forall x\in\mathbb{R}^d,t\geq0.
	\end{equation}
	
	Then for simplicity we assume the total density $n$ is expanding as a characteristic function for all time and we assume $\Omega(t)$ defined as in \eqref{eq:c-compact} satisfies:
	\begin{equation}
	\Omega(t)=\{x:n(x,t)>0\}=\{x:n(x,t)=1\}=\{x:p(x,t)>0\}
	\end{equation}
	
	Then from the equation of pressure $p$ \eqref{dev-con-relation}, we have
	\begin{equation*}
	\begin{cases}
	-\Delta p=\mu G(c)+(1-\mu)(G(c)-D),\ x\in\Omega(t),\\
	p=0,\ x\in\p \Omega(t),
	\end{cases}\\
	\end{equation*}where in the first equation we substitutes the definition of $G_1,G_2$ in \eqref{eq:g1g2}.
	
	Since $n\equiv1$ in $\Omega(t)$, we have $n_1=n\mu=\mu,n_2=n(1-\mu)=1-\mu$. Thus, from \eqref{eq:c-compact} we have
	\begin{equation*}
	\begin{cases}
	\epsilon\frac{\p c}{\p t}-\Delta c+\psi(c)=(1-\mu)a,\quad x\in\Omega(t),\\
	c=c_B,\quad x\in\p\Omega(t).
	\end{cases}
	\end{equation*}For simplicity, in the free boundary model we focus on the quasi-static case: $\epsilon=0$. We formally think the evolution of density fraction $\mu$ still satisfies (\ref{eq:model-mu2}) for $x\in\Omega(t)$ after taking the limit. 
	
	The last but subtlest component of the model is the evolution of $\Omega(t)$. We assume that the evolution of $\Omega(t)$ is determined by the motion of its boundary $\p\Omega(t)$ and that the moving speed of $\p\Omega(t)$ is still determined by Darcy's law. The moving speed in the normal direction satisfies:
	\begin{equation*}
	V_n=-\nabla p \cdot \textbf{n},\quad x\in\p\Omega(t).
	\end{equation*}where $\textbf{n}$ is the outer normal direction of $\Omega(t)$. 

	\subsubsection{Model and assumptions}
	To summarize, the free boundary model consists of 4 components: the time-dependent domain $\Omega(t)$, the pressure $p$ which determines the moving speed of $\p\Omega(t)$, nutrients $c$ and the density fraction(of normal cells) $\mu$ for two-cells dynamics, which corresponds to $\frac{n_1}{n_1+n_2}$ in the compressible model. If $\mu(x,t)=1$, it means locally all cells are normal cells. If $\mu(x,t)=0$, it means locally all cells are autophagic cells.
	
	Initially we shall specify a domain $\Omega(0)$ and an initial value for $\mu$, $\mu(x,0)=\mu_0(x)$ which is a function in $\Omega(0)$. Since $\mu_0$ is the density fraction, it should take values in $[0,1]$, i.e.,
	\begin{equation}\label{eq:free-init}
\mu(x,0)=\mu_0(x),\quad	\mu_0(x)\in[0,1],\quad x\in\Omega(0).
	\end{equation} 
	
	As time evolves, the density fraction $\mu$ satisfies a reaction-convection equation on a time-dependent domain while pressure $p$ and nutrient $c$ solves two elliptic equations for each time $t\geq0$.
		\begin{align}
		\label{eq:model-mu2}
		\frac{\partial \mu}{\partial t}-(\nabla \mu)\cdot(\nabla p)&=-\mu K_1(c)+(1-\mu)K_2(c)+D\mu(1-\mu),\  x\in\Omega(t).\\
		&\begin{cases}\label{eq:model-p}
		-\Delta p=\mu G(c)+(1-\mu)(G(c)-D),\ x\in\Omega(t),\\
		p=0,\ x\in\p \Omega(t).
		\end{cases}\\
		&\begin{cases}\label{eq:model-c}
		-\Delta c+\psi(c)=(1-\mu)a  ,\ x\in \Omega(t),\\
		c=c_B  ,\ x\in\p\Omega(t).
		\end{cases}
		\end{align} 
		
	The evolution of $\Omega(t)$ is determined by the motion of its boundary $\p\Omega(t)$. Same as the compressible model, the moving speed of boundary is governed by the negative gradient of pressure (Darcy's law). Precisely, the moving speed in normal direction is 		
	\begin{align}
	\label{eq:model-v}
	V_n=-\nabla p \cdot \textbf{n},\quad x\in\p\Omega(t).
	\end{align}
	
As we will show in Proposition \ref{prop:tumor-grow} in Section \ref{sc:math-proof}, under the Assumption \ref{as:non-negative} below, $V_n=-\nabla p \cdot \textbf{n}>0$ thus the characteristic is going outwards we do not need a boundary condition for the hyperbolic equation \eqref{eq:model-mu2}. 

Last but not least, we impose the following assumption, which is inherited from \eqref{dev:as2} in the formal derivation.
\begin{assumption}\label{as:non-negative}For $T>0$, the weighted net growth rate $\mu G+(1-\mu)(G-D)$ is non-negative, i.e.,
	\begin{equation}
\mu(x,t) G(c(x,t))+(1-\mu(x,t))(G(c(x,t))-D)\geq0,
	\end{equation} for all $x\in\Omega(t),0\leq t<T$. And for every $t\geq0$ there exists at least one point $x\in\Omega(t)$ such that the strict inequality holds.
\end{assumption}
Assumption \ref{as:non-negative} allows the net growth rate of autophagic cells $G(c)-D$ to be negative, but assumes the nutrient is sufficient enough such that the weighted net growth rate is non-negative at every point in tumor. This is important for the connection between the free boundary system \eqref{eq:free-init}$\sim$\eqref{eq:model-v}, which is already closed itself, and the incompressible limit of the compressible model \eqref{model:PDE-systems}. 

Note that if Assumption \ref{as:non-negative} is true, then by maximum principle \eqref{eq:model-p} yields that $p(x,t)>0$, for all $x\in\Omega(t)$. This rules out the formation of the necrotic core inside $\Omega(t)$, the case when the total density $n$ decays to be less than $1$ and pressure $p$ becomes zero (see section 7.3 in \cite{perthame2015some}). Biologically, necrotic core means dead cells form a core in the tumor center due to the lack of nutrients \cite{cui2001analysis,perthame2014traveling}. In the framework of incompressible limit of cell density model, the necrotic core is characterized as an obstacle problem in pressure, as shown in \cite{guillen2020heleshaw}. We will investigate in this direction in the future.

Besides, we assume the following regularity of solutions for the analysis in next section.
\begin{assumption}\label{as:strong-simple-con}
	\begin{enumerate}For $T>0$, we have
		\item $\Omega(t)$ is a bounded and simple-connected domain with smooth boundary for all $t\in(0,T)$. And the moving speed of boundary satisfies (\ref{eq:model-v}), that is, there exists a parameterization of the boundary $\Omega(t)$: $x(t,\alpha),\ \alpha\in[0,1],$ which satisfies
		\begin{equation}
		\frac{d}{dt}x(t,\alpha)=-\nabla p(x(t,\alpha),t),\quad x(t,0)=x(t,1).
		\end{equation}
		for all $t\in(0,T),\ \alpha\in[0,1]$. 
		\item Let $Q_T=\{(x,t),x\in\overline{\Omega(t)},\ t\in[0,T)\}$. Then $p,c,\mu\in C^2(Q_T)$. And  \eqref{eq:model-mu2}$\sim$\eqref{eq:model-c} are satisfied in classical sense. 
	\end{enumerate}
\end{assumption}The global and explicit formulation of moving speed of boundary is convenient for the flow map formulation of $\Omega(t)$ (Proposition \ref{prop:refor}) in next section.
\hl{
\paragraph{Remark on well-posedness of the free boundary model} Although the incompressible limit has been established in \cite{liu2021existence}, it remains difficult show that the limit solution is exactly the solution of the free boundary problem (like one species case \cite{ahelemellet2017hele,kim2018porous}). Nevertheless, one may try to prove the well-posedness of the free boundary model directly. In this direction, the most relevant work in literature might be \cite{general2003} which also treats a multi-species free boundary model coupled with nutrients. However, there is a crucial difference: in \cite{general2003} (and related literature) the boundary value of pressure $p$ is proportional to the mean curvature $\kappa$ of $\p\Omega(t)$:
\begin{equation}
p=\gamma_0\kappa,\quad x\in\p\Omega(t).
\end{equation}The coefficient $\gamma_0>0$ denotes the surface tension. While in our free boundary model, derived from a cell density model, the pressure is zero on the boundary, which can be seen as the case $\gamma_0=0$. If we add the surface tension to our model, then the approach in \cite{general2003} can be adapted quite directly to obtain local well-posedness. But the mean curvature as an elliptic operator has a regularizing effect and somehow makes the problem parabolic. Therefore it may be difficult to directly use the method in \cite{general2003}, especially extend the Theorem 2.1 in \cite{general2003} to our case $\gamma_0=0$. A natural approach is to establish the well-posedness for $\gamma_0>0$ then let $\gamma_0\rightarrow0$, which is so-called vanishing surface tension limit in literature \cite{hwang2016vanishing,yi1997quasi}. If one can combine the treatment for multi-species model \cite{general2003} and the vanishing surface tension limit in \cite{yi1997quasi}, one may obtain the local well-posedness of our model. }

\hl{
An alternative approach may be to learn from direct proofs of well-posedness for the Hele-Shaw problem when $\gamma=0$. Early results focus on two dimension \cite{vinogradov1948problem,gustafsson1984differential,reissig1993simplified}. \cite{escher1997classical} first prove the existence of Hele-Shaw problem when the initial boundary $\p\Omega(t)$ is $C^{2+\alpha}(\alpha\in(0,1))$. This may be the closest paper to get a $\gamma_0=0$ version of Thm 2.1 in \cite{general2003} since \cite{general2003} also uses a H\"{o}lder setting. We remark that there are also viscosity solution approach \cite{kim2003uniqueness} and weak formulation \cite{elliott1981variational} for the well-posedness of Hele-Shaw problem without surface tension.} 

\hl{
However, these are beyond the scope of this paper, since rigorously proving the well-posedness of the free boundary model is not the focus of the current work, and is worth separate exploration.}

\hl{
Last but not least, we remark that in our case surface tension $\gamma_0=0$, Assumption \ref{as:non-negative} shall be crucial for the well-posedness. As aforementioned, this assumption ensure $\Omega(t)$ is expanding. When $\Omega(t)$ is retreating the Hele-Shaw problem without surface tension (i.e. $\gamma_0=0$) is known to be ill-posed \cite{gustafsson2006conformal}.
}
\section{Analysis on the free boundary model and the well mixed limit}\label{sc:math-proof}

		This section is devoted to analyzing the free boundary model \eqref{eq:free-init}$\sim$\eqref{eq:model-v}. First we give basic properties which shed some light on the structure of solutions, then we study the asymptotic behavior of the density fraction $\mu$. In a special case, we show that the density fraction $\mu$ converges to a spatial homogeneous state within the tumor region $\Omega(t)$ as $t$ goes to infinity. Moreover, the convergence is exponential. We call this phenomenon the  ``well-mixed'' limit. 
	\subsection{Well-posedness, characteristic structure and boundness of $\mu$}
	In this subsection, we give some basic properties for the free boundary model. In particular we prove the characteristic structure of the density fraction $\mu$ which is useful in analyzing its asymptotic behavior.
	
	First, we have the following proposition on the moving speed of $\p\Omega(t)$.
	\begin{proposition}\label{prop:tumor-grow}
For $T>0$, suppose $(\Omega,p,c,\mu)$ is a solution of the free boundary system \eqref{eq:free-init}$\sim$\eqref{eq:model-v} which satisfies Assumption \ref{as:non-negative} and \ref{as:strong-simple-con} on the time interval $(0,T)$.
		Then the moving speed of the boundary $\p\Omega(t)$ in normal direction is positive, that is,
		\begin{equation*}
		V_n=-\nabla p\cdot \mathbf{n}>0,
		\end{equation*}for every $x\in\p\Omega(t)$ and $t\in(0,T)$.
	\end{proposition}
	\begin{proof}From Assumption \ref{as:non-negative} we have $\Delta p= -(\mu G+(1-\mu)(G-D))\leq0 $ for all $x\in\Omega(t)$, thus $p$ achieves its minimum in $\overline{\Omega(t)}$ at $\p\Omega(t)$. Thanks to the Dirichlet zero boundary condition, every point on $\partial\Omega(t)$ is a minimum point. Thus by Hopf's lemma we have $\frac{\partial p}{\partial n}<0$ for all $x\in\p\Omega(t)$, unless $p(\cdot,t)\equiv
		0$ in $\Omega(t)$. The latter case implies $\Delta p\equiv0$ which contradicts the Assumption \ref{as:non-negative}.
	\end{proof}

Proposition \ref{prop:tumor-grow} ensures that we do not need a boundary condition for the density fraction $\mu$. To give a clear statement of its useful characteristic structure, we first give a reformulation of $\Omega(t)$ through the following flow map:
	\begin{equation}\label{eq:flowmap}
	\begin{cases}
	\frac{\mathrm{d}x}{\mathrm{d}t}(t,y)=-\nabla p(x(t,y),t).\\
	x(0,y)=y\in\overline{\Omega(0)}.
	\end{cases}
	\end{equation}Well-posedness of this flow map is stated in the following proposition.
	\begin{proposition}\label{prop:refor}
	For $T>0$, suppose $(\Omega,p,c,\mu)$ is a solution of the free boundary system \eqref{eq:free-init}$\sim$\eqref{eq:model-v} which satisfies Assumption \ref{as:non-negative} and \ref{as:strong-simple-con} on the time interval $(0,T)$. Then the followings are true.
		\begin{enumerate}
			\item The solution of the flow map (\ref{eq:flowmap}) $x(t,y)$ exists and is unique for all initial $y\in \overline{\Omega(0)}, t\in[0,T)$.
			\item Set $X_{t,0}$ be the map defined on $\overline{\Omega(0)}$ such that $X_{t,0}(y)=x(t,y)$ for $t\geq0$. Then $X_{t,0}$ is a homeomorphism from $\overline{\Omega(0)}$ to $\overline{\Omega(t)}$. In particular for all initial $x\in\overline{\Omega(t)}$ there exists a unique $y\in\overline{\Omega(0)}$ such that $x=X_{t,0}(y)=x(t,y)$.
		\end{enumerate}
	\end{proposition}

	Note well-posedness of \eqref{eq:flowmap} for $x$ at $\p\Omega(t)$ is implied by Assumption \ref{as:strong-simple-con}. The proof of Proposition \ref{prop:refor} is straightforward and we postpone it to Appendix \ref{app:proof}. With this flow map formulation of $\Omega(t)$, we could state the characteristics structure of $\mu$. Consider the following characteristic ODE system for the quasi-linear equation \eqref{eq:model-mu2}:
	\begin{equation}\label{eq:char-mu}\begin{cases}
	&\frac{dx}{dt}(t,y)=-(\nabla p)(x(t,y),t),\\
	&\frac{dz}{dt}(t,y)=-z K_1+(1-z)K_2+Dz(1-z),\\
	&x(0,y)=y\in\overline{\Omega(0)},\\
	&z(0,y)=\mu_0(y).
	\end{cases}
	\end{equation}Here the transition rates $K_1,K_2$ in the second equation depend on the nutrient at point $x(t,y)$ and time $t$, i.e., $K_i=K_i(c(x(t,y),t))$, $i=1,2$. We use the shorthand notation to save space. In \eqref{eq:char-mu} $z$ can be seen as the Lagrangian representation of the density fraction $\mu$. This fact is stated in the following proposition.
	\begin{proposition}\label{prop:char}
				For $T>0$, suppose $(\Omega,p,c,\mu)$ is a solution of the free boundary system \eqref{eq:free-init}$\sim$\eqref{eq:model-v} which satisfies Assumption \ref{as:non-negative} and \ref{as:strong-simple-con} on the time interval $(0,T)$.
		Then for all $x\in\overline{\Omega(t)},t\in(0,T)$ there exists a $y\in\overline{\Omega(0)}$ such that $x=x(t,y)$ and we have $\mu(x,t)=z(t,y)$. Here $x,z$ are the characteristic ODE as in \eqref{eq:char-mu}. 
	\end{proposition}
	\begin{proof}
		Proposition \ref{prop:refor} ensures that for all $x\in\Omega(t)$ there exists a unique $y\in\Omega(0)$ such that $x=x(t,y)$. Then the following computation is standard:
		\begin{align*}
		\frac{d}{ds}\mu(x(s,y),s)&=\frac{\partial \mu}{\partial t}-(\nabla \mu)\cdot(\nabla p)\\&=-\mu K_1+(1-\mu)K_2+D\mu(1-\mu).
		\end{align*}
		Thus $\mu(x(\cdot,y),\cdot)$ shares the same ODE and the initial value with $z(\cdot,y)$ in \eqref{eq:char-mu}, and one obtains $\mu(x,t)=z(t,y)$.
	\end{proof}
Proposition \ref{prop:char} in helpful in the proof of the uniform convergence of the ``well-mixed'' limit in next subsection.
 With this characteristic structure, now it is easy give a point-wise bound for $\mu$ as follows.
	\begin{corollary}\label{cor:01}
				For $T>0$, suppose $(\Omega,p,c,\mu)$ is a solution of the free boundary system \eqref{eq:free-init}$\sim$ \eqref{eq:model-v} which satisfies Assumption \ref{as:non-negative} and \ref{as:strong-simple-con} on the time interval $(0,T)$. And the parameters satisfy \eqref{ap-2GK}.
		  Then the density fraction $\mu(x,t)\in[0,1]$, for all $x\in\overline{\Omega(t)}$ and  $t\in(0,T)$.
	\end{corollary}
	\begin{proof}
		Let us look at $z(t,y)$, the Lagrangian representation of $\mu$  in \eqref{eq:char-mu}. By Proposition \ref{prop:char} it suffices to show that $z(t,y)\in[0,1]$ for all $y\in\overline{\Omega(0)}$ and $t\in(0,T)$. Recall the ODE for $z$ \eqref{eq:char-mu}:
		\begin{equation*}
		\frac{dz}{dt}=-zK_1+(1-z)K_2+Dz(1-z).
		\end{equation*}Notice that initially $z(0,y)=\mu(y,0)=\mu_0(y)\in[0,1]$ by \eqref{eq:free-init}. Then the result follows from that $K_1,K_2>0$ in \eqref{ap-2GK} and the comparison principle.
	\end{proof}
	\subsection{Asymptotic behavior for $\mu$: the well-mixed limit}\label{subsc:asymp-mu}
	This subsection is devoted to analyzing the asymptotic behavior of density fraction $\mu$ as time $t$ goes to infinity. For simplicity of analysis, we make the following assumption.
	\begin{assumption}\label{as:constants} The transition rates $K_1,K_2$ are constants.
	\end{assumption}

With Assumption \ref{as:constants}, we can show that the density fraction $\mu$ will converge to a spatial-homogeneous steady state $\mu\equiv\mu^*,\ x\in\Omega(t)$, where $\mu^*\in[0,1]$ is a constant.  And we characterize it from two aspects: a uniform convergence result (Theorem \ref{thm:uni-con}) and convergence in ``$L^{2n}$ norm''($n\in\mathbb{N}^+$) under a family of conditions (Theorem \ref{thm:L2n}).
	
	First we analyze the reaction term for $\mu$. Let $f(\mu)$ denotes the right hand side in \eqref{eq:model-mu2}, we have
	\begin{equation}\label{eq:rec-ode}
	f(\mu)\coloneqq-\mu K_1+(1-\mu)K_2+D\mu(1-\mu)=-D\mu ^2+(D-K_1-K_2)\mu+K_2.
	\end{equation}
	 We observe that $f$ is a quadratic function in $\mu$. Recall in our model Assumption \ref{ap-1Da} the extra death rate due to autophagy $D>0$. And in \eqref{ap-2GK} the transition rates $K_1,K_2>0$. Thus we observe that $f(0)=K_2>0,f(1)=-K_1<0$ and $-\frac{K_2}{D}<0$. This implies that the two roots $\nu^*<\mu^*$ of $f$ satisfies that $\nu^*<0,\ \mu^*\in(0,1).$ Precisely we have
	 \begin{equation}\label{eq:roots}
	 \nu^*=\frac{D-K_1-K_2-\sqrt{E}}{2D},\quad\mu^*=\frac{D-K_1-K_2+\sqrt{E}}{2D},
	 \end{equation}where $E=D^2+(K_1+K_2)^2-2DK_1+2DK_2$.

	With this information we could characterize the asymptotic behavior of the reaction ODE $\frac{dz}{dt}=f(z)$, which we write as the following lemma.
	
	\begin{lemma}\label{lemma:ode}Suppose \eqref{ap-2GK},\eqref{ap-1Da} and Assumption \ref{as:constants} are satisfied. Let $f(\mu)$ be the rate function as in \eqref{eq:rec-ode} and $\nu^*<\mu^*$ be two roots of $f$ as in \eqref{eq:roots}. Consider the following ODE,
		\begin{equation*}
		\begin{cases}
		\frac{dz}{dt}(t)=f(z(t))=-D(z(t)-\nu^*)(z(t)-\mu^*),\ t\geq 0.\\
		z(0)=z_0.
		\end{cases}
		\end{equation*}If initial value $z_0>\nu^*$, then the solution $z(t)$ converges to $\mu^*$ exponentially as $t\rightarrow\infty$. Precisely, we have 
		\begin{equation}\label{eq:ode-wellmixed}
		|z(t)-\mu^*|\leq\frac{\max\{z_0,\mu^*\}-\nu^*}{z_0-\nu^*}e^{-D(\mu^*-\nu^*)t}|z_0-\mu^*|.
		\end{equation}
	\end{lemma}
\begin{proof}From the comparison principle we know that if $\mu^*\leq z_0$ then $\nu^*<\mu^*\leq z(t)$, for all $t\geq0$, and thus $\frac{dz}{dt}(t)=-D(z(t)-\nu^*)(z(t)-\mu^*)\leq0$, for all $t\geq0$, therefore the solution is decreasing, in particular $z(t)\leq z_0$, for all $t\geq0$. 
	
	Otherwise if $z_0\in(\nu^*,\mu^*)$ then from the comparison principle $\nu^* \leq z(t)\leq \mu^*$ for all $t\geq0$. 
	
	In both cases one has 
	\begin{equation*}
	\nu^*\leq z(t)\leq\max\{z_0,\mu^*\},\ \forall t\geq0.
	\end{equation*}
	On the other hand with elementary calculation we have
	\begin{equation*}
	{|z(t)-\mu^*|}=\frac{{|z(t)-\nu^*|}}{|z_0-\nu^*|}e^{-D(\mu^*-\nu^*)t}|z_0-\mu^*|.
	\end{equation*}
	Then the result follows from $0<z(t)-\nu^*\leq\max\{z_0,\mu^*\}-\nu^*.$
\end{proof}
Applying Lemma \ref{lemma:ode} along the characteristics of density fraction $\mu$ \eqref{eq:char-mu}, we obtain the following uniform convergence result:
	\begin{theorem}\label{thm:uni-con}
		Suppose $(\Omega,p,c,\mu)$ is a solution of the free boundary system \eqref{eq:free-init}$\sim$\eqref{eq:model-v} which satisfies Assumption \ref{as:non-negative} and \ref{as:strong-simple-con} on the time interval $(0,\infty)$. Suppose \eqref{ap-1Da}, \eqref{ap-2GK} and Assumption \ref{as:constants} are satisfied. Let $f$ be the rate function as in \eqref{eq:rec-ode} and $\nu^*<\mu^*$ be two roots of $f$ as in \eqref{eq:roots}. 
		
		Then we have $\mu(x,t)\rightarrow \mu^* $ uniformly for $x\in\overline{\Omega(t)}$ as $t\rightarrow \infty$.
	Here ``uniformly'' means $||\mu(\cdot,t)-\mu^*||_{C(\overline{\Omega(t)})}\rightarrow0$ as $t\rightarrow \infty$. In addition, the convergence is exponentially fast in the sense that \begin{equation}\label{eq:uni-exp}
	||\mu(\cdot,t)-\mu^*||_{C(\overline{\Omega(t)})} \leq A e^{-D(\mu^*-\nu^*)t}||\mu(\cdot,0)-\mu^*||_{C(\overline{\Omega(0)})},\end{equation} where $A>0$ is some constant.
	\end{theorem}
	\begin{proof}
		By Proposition \ref{prop:char}, for all $x$ in $\overline{\Omega(t)}$, despite $\Omega(t)$ is expanding, we can always trace back along characteristics to find an initial position $y\in\overline{\Omega(0)}$ such that 
		\begin{equation*}
		\mu(x,t)=z(t,y),\ x=x(t,y),
		\end{equation*}where $x(t,y)$ is the characteristics of $\mu$ and $z(t,y)$ is the Lagrangian representation for $\mu$ as in \eqref{eq:char-mu}.
		Since with initial condition \eqref{eq:free-init}, $z(0,y)=\mu(y,0)=\mu_0(y)\in[0,1]\subset(\nu^*,+\infty)$, by Lemma \ref{lemma:ode} we have 
		\begin{align*}
		|z(t,y)-\mu^*|&\leq\frac{\max\{z(0,y),\mu^*\}-\nu^*}{z(0,y)-\nu^*}e^{-D(\mu^*-\nu^*)t}|z(0,y)-\mu^*|\\
		&\leq \max_{z_0\in[0,1]}\{\frac{\max\{z_0,\mu^*\}-\nu^*}{z_0-\nu^*}\}e^{-D(\mu^*-\nu^*)t}|z(0,y)-\mu^*|,
		\end{align*}for all $y\in\overline{\Omega(0)}.$ With straightforward computation, one has
		\begin{equation*}
\frac{\max\{z_0,\mu^*\}-\nu^*}{z_0-\nu^*}=\frac{\max\{z_0-\nu^*,\mu^*-\nu^*\}}{z_0-\nu^*}=\max\{1,\frac{\mu^*-\nu^*}{z_0-\nu^*}\}\leq\frac{\mu^*-\nu^*}{-\nu^*},
		\end{equation*} when $z_0\in[0,1]$. Then we conclude
		\begin{equation*}
		|z(t,y)-\mu^*|\leq \frac{\mu^*-\nu^*}{-\nu^*} e^{-D(\mu^*-\nu^*)t}|z(0,y)-\mu^*|.
		\end{equation*} Taking maximum for $y\in\overline{\Omega(0)}$ at both sides, we complete the proof. And the coefficient $A>0$ in \eqref{eq:uni-exp} can be taken as $\frac{\mu^*-\nu^*}{-\nu^*}.$
	\end{proof}
	Note that the convergence can also be observed directly through the equation of $\mu$ (\ref{eq:model-mu2}), which is equivalent to
	\begin{equation*}
	\partial_t(\ln\frac{|\mu-\mu^*|}{|\mu-\nu^*|})-\nabla p\cdot(\nabla \ln\frac{|\mu-\mu^*|}{|\mu-\nu^*|})=-D{(\mu^*-\nu^*)}.
	\end{equation*}
	
\begin{figure}
		\begin{subfigure}{1.0\textwidth}
			\centering
		\includegraphics[width=1.0\textwidth]{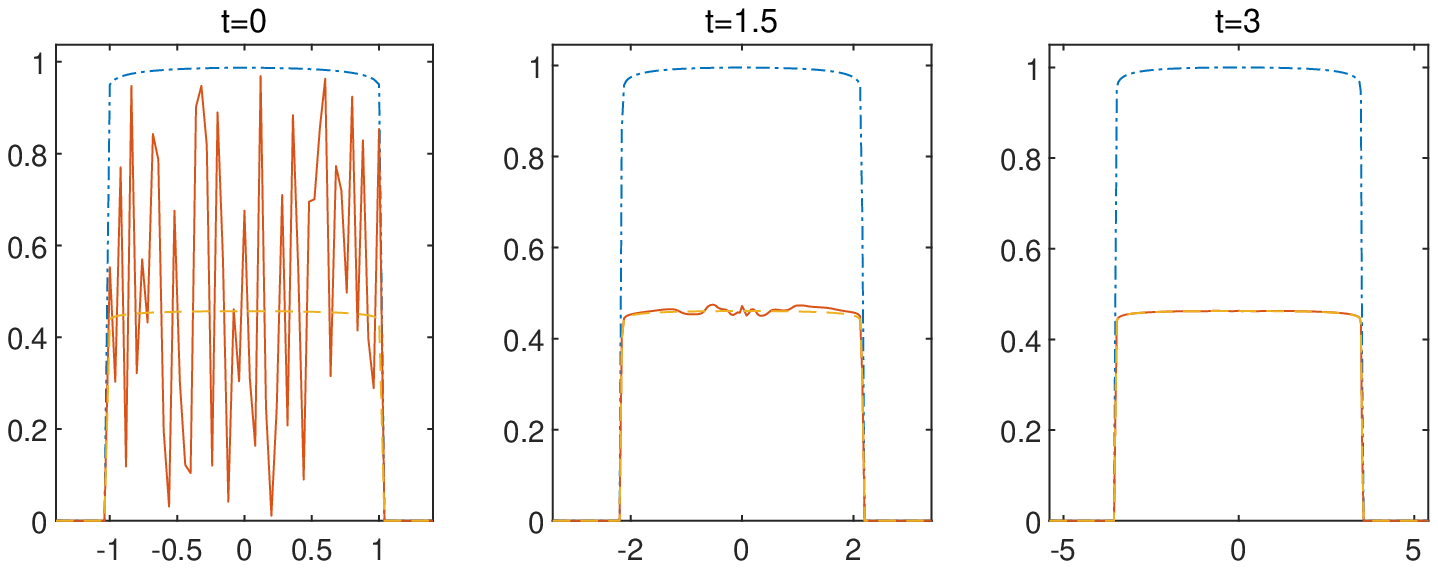}
		\caption{$\gamma=80$}
		\label{figure2a}
	\end{subfigure}
	\begin{subfigure}{1.0\textwidth}
		\centering
		\includegraphics[width=1.0\textwidth]{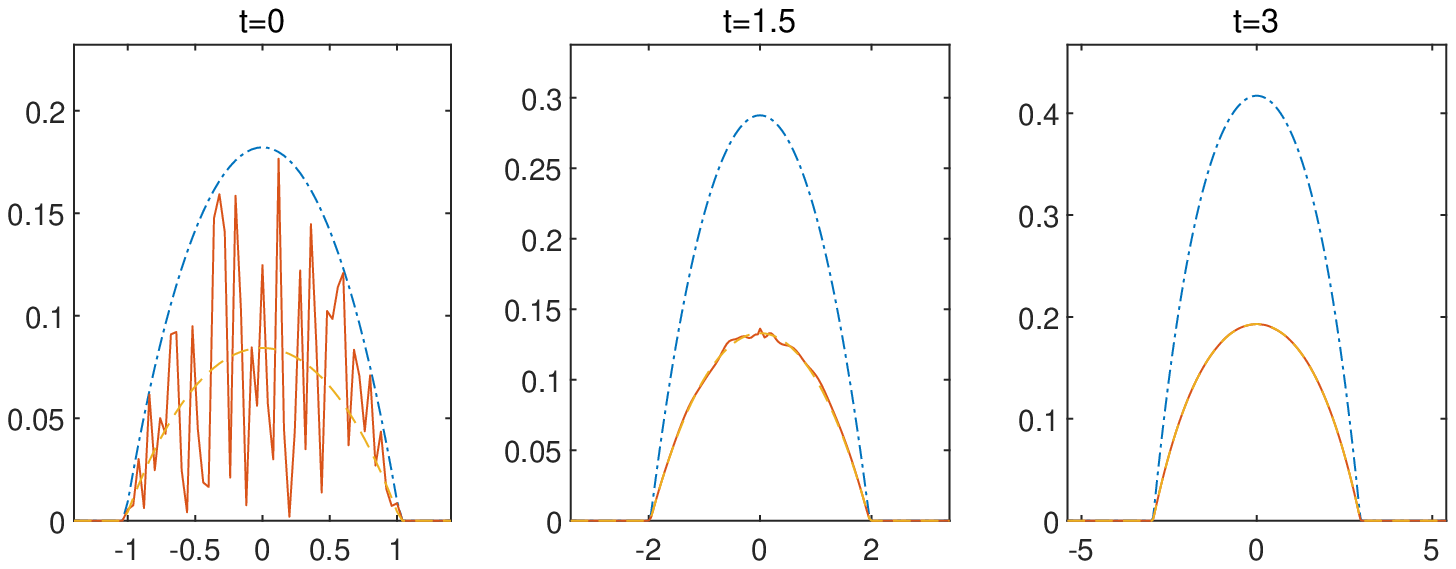}
		\caption{\hl{$\gamma=2$}}
		\label{figure2b}
	\end{subfigure}
	\caption{Densities of cells at time $t=0,1.5,3$ with a heterogeneous initial density fraction. Blue dash-dotted line: the total density $n$. Red solid line: the density of autophagic cells $n_2$. Orange dashed line: theoretical value for density of autophagic cells in the ``well-mixed'' limit. Parameters: $G(c)=gc,\psi(c)=c,K_1=K_2=1,g=1,a=0.4,D=0.3$.}
	\label{fig:s3unicon}
\end{figure}

	The uniform convergence result is illustrated in Figure \ref{figure2a}, where we numerically simulate the compressible model with large $\gamma$. We plot the densities of total cells $n$ and autophagic cells $n_2=(1-\mu)n$. We also plot the theoretical value for the density of autophagic cells in the ``well-mixed'' limit $(1-\mu^*)n$ for comparison. Initially the value of $\mu$ is highly spatially heterogeneous, and as time evolves, the tumor expands while the two types of cells mix towards the equilibrium ratio $\mu^*$. In particular when $t=3$ autophagic cells and normal cells are effectively ``well-mixed''.
	
		\hl{ This well-mixed limit is expected to hold also for the cell density model \eqref{model:PDE-systems}, since it also has the $\mu$ equation \eqref{eq:model-mu-Com}. Actually if the solution is smooth, then the above characteristic argument still works. The smoothness of solution may be ensured if initial data is smooth and bounded from below. For general case a regularizing argument may be needed, which is beyond the scope of this paper. Nevertheless we simulate the case $\gamma=2$ for illustration in Figure \ref{figure2b}.}

	Besides this uniform convergence result, we consider the ``$L^q$ estimate'' on $||\mu-\mu^*||_{L^q(\Omega(t))}(1\leq q<+\infty)$. If the domain is fixed, then the decay of $\mu-\mu^*$ in $L^q$ norm is trivial since the $L^\infty$ norm can control the $L^q$ norm. The difference is here we are facing an expanding domain. We first give an intuitive analysis. From Theorem \ref{thm:uni-con} the ``pointwise'' convergence of $|\mu-\mu^*|$ is exponential with the rate $r_1:=D(\mu^*-\nu^*)$. Thus, intuitively one has
	\begin{equation}\label{eq:expdecay1}
	\int_{\Omega(t)}|\mu-\mu^*|^{q}dx\approx C|\Omega(t)|\exp\{-qr_1t\}.
	\end{equation} As we will see in the next section (see discussion about \eqref{v1},\eqref{eq:exp-growth}), when the nutrient provided by autophagy is relatively sufficient, $|\Omega(t)|$ may also grow exponentially, i.e.,
	 \begin{equation}\label{eq:expgrow2}
	 |\Omega(t)|\approx|\Omega(0)|\exp\{r_2t\},
	 \end{equation}for some $r_2>0$. Combine \eqref{eq:expgrow2} with \eqref{eq:expdecay1}, formally we have
	 \begin{equation}\label{eq:critical}
	 \int_{\Omega(t)}(\mu-\mu^*)^{q}dx\approx C|\Omega(0)|\exp\{(r_2-qr_1)t\}.
	 \end{equation} Thus we are facing a competition of exponential growth and exponential decay. Intuitively the $L^q$ norm ($q\geq1$) $||\mu-\mu^*||_{L^q(\Omega(t))}$ will exponentially decay if and only if $q>\frac{r_2}{r_1}$. In the following, we do not investigate \eqref{eq:critical} in details but give a family of sufficient conditions for the decay in ``$L^{2n}$ norm''($n\in\mathbb{N}^+$) to be held. For preparation, first we give an a priori estimate on the concentration of nutrients $c$.
	 \begin{lemma}\label{lemma:c}
For $T>0$, suppose $(\Omega,p,c,\mu)$ is a solution of the free boundary system \eqref{eq:free-init}$\sim$\eqref{eq:model-v} which satisfies Assumption \ref{as:non-negative} and \ref{as:strong-simple-con} on the time interval $(0,T),T>0$. And suppose \eqref{ap-1Da}$\sim$\eqref{ap-3psi} are true. Then we have
$$c(x,t)\leq\max\{c_B,c_0\},$$for all $x\in\overline{\Omega(t)}$, $t\in(0,T)$. Recall $c_B$ is the boundary value of $c$ and $c_0$ is the threshold concentration such that $\psi(c_0)=a$ as in \eqref{ap-3psi}. 
	 \end{lemma}
 The proof is straightforward from the maximum principle.
\begin{proof}
For $t\in(0,T)$, suppose $c(\cdot,t)$ reaches it maximum in $\overline{\Omega(t)}$ at $x^*\in\overline{\Omega(t)}$. It suffices to prove $c^*:=c(x^*,t)\leq\max\{c_0,c_B\}$. If $x^*\in\p\Omega(t)$, then $c^*=c_B$. It is done.

Otherwise the maximum point $x^*$ lies in $\Omega(t)$, and one obtains $\Delta c(x^*,t)\leq0$. Thus from the equation \eqref{eq:model-c} for $c$, one gets
\begin{equation}\label{eq:proof-c01}
\psi(c^*)\leq(1-\mu)a\leq a=\psi(c_0).
\end{equation} The second inequality uses that $\mu(x,t)\in[0,1]$ by Corollary \ref{cor:01}. And the last equality is the definition of $c_0$ \eqref{ap-3psi}. Thus we obtain $c^*\leq c_0$ since $\psi'(c)>0$ as in \eqref{ap-3psi}.
\end{proof}
 Thanks to Lemma \ref{lemma:c}, we could control the concentration of nutrients $c$. Now we give a family of sufficient conditions for $||\mu-\mu^*||_{L^{2n}(\Omega(t))}(n\in\mathbb{N}^+)$ to exponentially decay:
	\begin{theorem}\label{thm:L2n}
	Suppose $(\Omega,p,c,\mu)$ is a solution of the free boundary system \eqref{eq:free-init}$\sim$\eqref{eq:model-v} which satisfies Assumption \ref{as:non-negative} and \ref{as:strong-simple-con} on the time interval $(0,\infty)$. And suppose \eqref{ap-1Da}$\sim$\eqref{ap-3psi} and Assumption \ref{as:constants} are true. Let $f$ be the rate function as in \eqref{eq:rec-ode} and $\nu^*<\mu^*$ be two roots of $f$ as in \eqref{eq:roots}. Then for $n\in\mathbb{N}^+$, $||\mu(\cdot,t)-\mu^*||_{L^{2n}(\Omega(t))}$ will exponentially decays if  
	\begin{equation}\label{eq:condition-thm2}
	G(\max\{c_B,c_0\})-D<2nK_2.
	\end{equation}Precisely, we have
		\begin{equation*}
		||\mu(\cdot,t)-\mu^*||_{L^{2n}(\Omega(t))}\leq\exp\{-Ct\}||\mu(\cdot,0)-\mu^*||_{L^{2n}(\Omega(0))},
		\end{equation*}where 
		\begin{equation*}
		C=\frac{-\nu^*}{1-\nu^*}K_1+\frac{1}{2n}(2nK_2-G(\max\{c_0,c_B\})+D)>0.
		\end{equation*}
	\end{theorem}
	\begin{proof}
		Let $\phi(\mu)=(\mu-\mu^*)^{2n}$, then by Reynolds's transport theorem,
		\begin{equation*}
		\frac{d}{dt}\int_{\Omega(t)}\phi(\mu)dx=\int_{\Omega(t)}\phi'(\mu)\mu_tdx+\int_{\partial\Omega(t)}\phi(\mu)(-\frac{\partial p}{\partial n})dS.
		\end{equation*}
		Since by the divergence theorem
		\begin{equation*}
		\int_{\partial\Omega(t)}\phi(\mu)(-\frac{\partial p}{\partial n})dS=\int_{\Omega(t)}-\phi'(\mu)\nabla \mu\cdot \nabla p-\phi(\mu)\Delta pdx,
		\end{equation*}
		we have
		\begin{align}
		\frac{d}{dt}\int_{\Omega(t)}\phi(\mu)dx&=\int_{\Omega(t)}\phi'(\mu)(\mu_t-\nabla \mu\cdot \nabla p)dx+\int_{\Omega(t)}\phi(\mu)(-\Delta p)dx\notag\\
		&=\int_{\Omega(t)}(\phi'(\mu)f(\mu)-\phi(\mu)\Delta p)dx,\label{eq:l2n-p1}
		\end{align}
		where $f$ is as in \eqref{eq:rec-ode}.
		Notice that
		\begin{equation*}
		\phi'(\mu)f(\mu)=-2nD(\mu-\mu^*)^{2n}(\mu-\nu^*)=-2nD\phi(\mu)(\mu-\nu^*).
		\end{equation*}
		Plugging it in \eqref{eq:l2n-p1} and substituting the equation for $p$ \eqref{eq:model-p} one obtains
		\begin{align}
		\frac{d}{dt}\int_{\Omega(t)}\phi(\mu)dx&=\int_{\Omega(t)}\phi'(\mu)f(\mu)-\phi(\mu)\Delta pdx\notag\\
		&=\int_{\Omega(t)}-\phi(\mu)(2nD(\mu-\nu^*)-G(c)\mu-(G(c)-D)(1-\mu))dx\notag\\
		&=-\int_{\Omega(t)}\phi(\mu)(2nD(\mu-\nu^*)-D\mu-(G(c)-D))dx,\label{eq:l2n-p2}
		\end{align}By Lemma \ref{lemma:c} we have $c(x,t)\leq\max\{c_0,c_B\}$, thus $G(c)-D\leq G(\max\{c_0,c_B\})-D=:G_2^*$. Applying this inequality to \eqref{eq:l2n-p2}, we have
		\begin{equation}\label{eq:l2n-p3}
		\frac{d}{dt}\int_{\Omega(t)}\phi(\mu)dx\leq-\int_{\Omega(t)}\phi(\mu)(2nD(\mu-\nu^*)-D\mu-G_2^*)dx
		\end{equation}
		Thanks to Gronwall's inequality, it suffices to show that $$\Phi(\mu)\triangleq2nD(\mu-\nu^*)-D\mu-G_2^*$$ has a positive lower bound.
		
		Note that $\Phi(\mu)$ is a linear function of $\mu\in[0,1]$ with slope $\Phi'(\mu)=(2n-1)D>0$. Thus it is equivalent to show that
		\begin{equation*}\label{pf:eq-l2n}
		\Phi(0)=-G_2^*-2nD\nu^*>0.
		\end{equation*}
		Recall that $\nu^*$ is a negative root of $f$:
		\begin{equation*}
		D\nu^*(1-\nu^*)-K_1\nu^*+K_2(1-\nu^*)=0.
		\end{equation*}
		Thus, multiply $\Phi(0)$ by $(1-\nu^*)>0$, we conclude that
		\begin{align*}
		\Phi(0)(1-\nu^*)&=-2nK_1\nu^*+2nK_2(1-\nu^*)-(1-\nu^*)G_2^*\\
		&=(1-\nu^*)(2nK_2-G_2^*)-2nK_1\nu^*>0.
		\end{align*}
		Then we obtain a positive lower bound for $\Phi$
		\begin{equation*}
		\Phi(\mu)\geq\Phi(0)\geq(2nK_2-G_2^*)+\frac{-\nu^*}{1-\nu^*}2nK_1>0,\quad \forall\mu\in[0,1].
		\end{equation*}Apply this to \eqref{eq:l2n-p3}, we obtain
		\begin{equation*}
	\frac{d}{dt}\int_{\Omega(t)}\phi(\mu)dx\leq-C\int_{\Omega(t)}\phi(\mu)dx,
		\end{equation*}where $C=(2nK_2-G_2^*)+\frac{-\nu^*}{1-\nu^*}2nK_1$. Then by Gronwall's inequality the proof is completed.
	\end{proof}

\begin{figure}[H]
		\centering
		\begin{subfigure}{0.49\linewidth}
			\centering
		\includegraphics[width=0.9\linewidth]{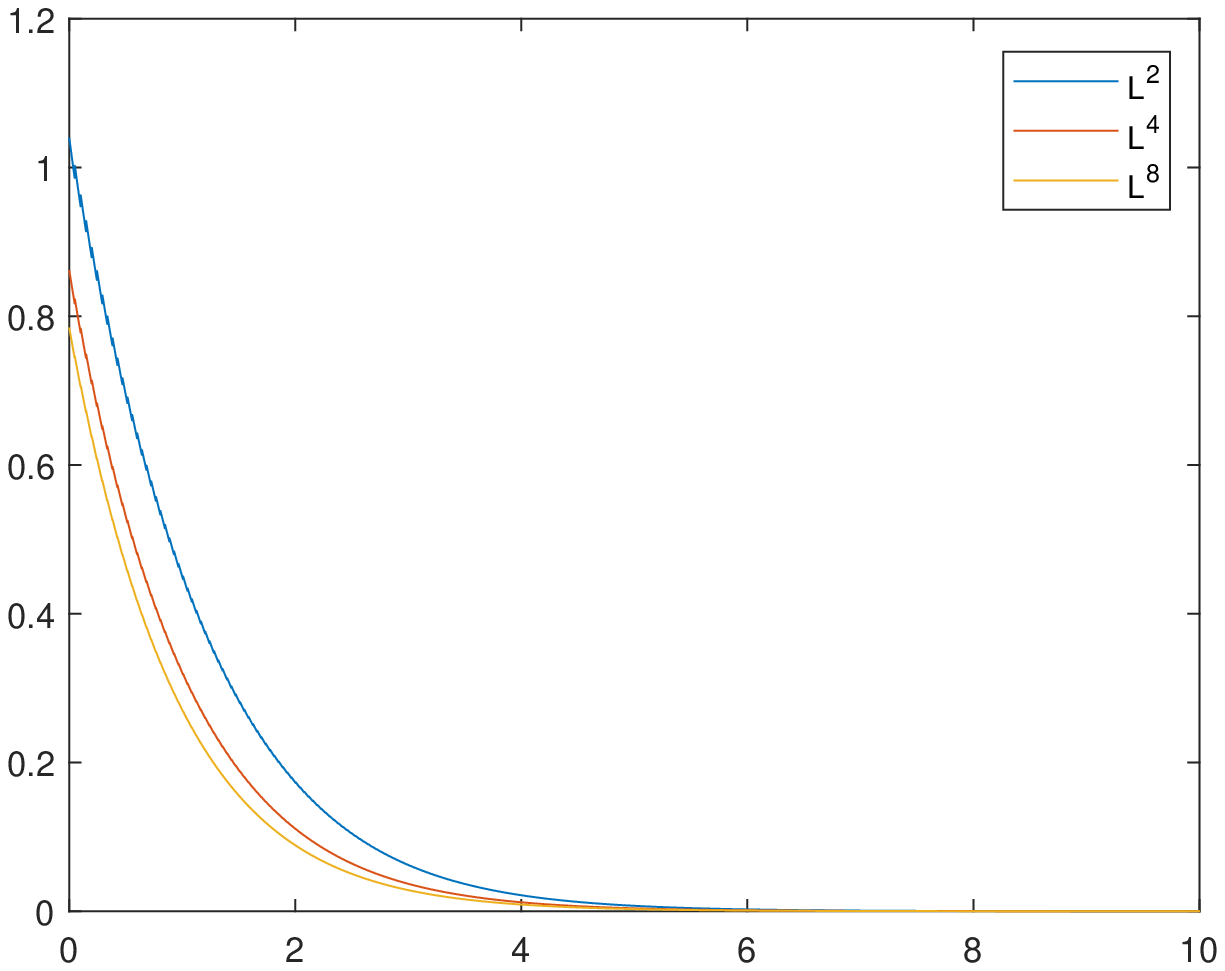}
\caption{$D=0.1,K_1=0.1,K_2=1$}
		\label{fig:l2n1}
		\end{subfigure}
\begin{subfigure}{0.49\linewidth}
	\centering
		\includegraphics[width=0.9\linewidth]{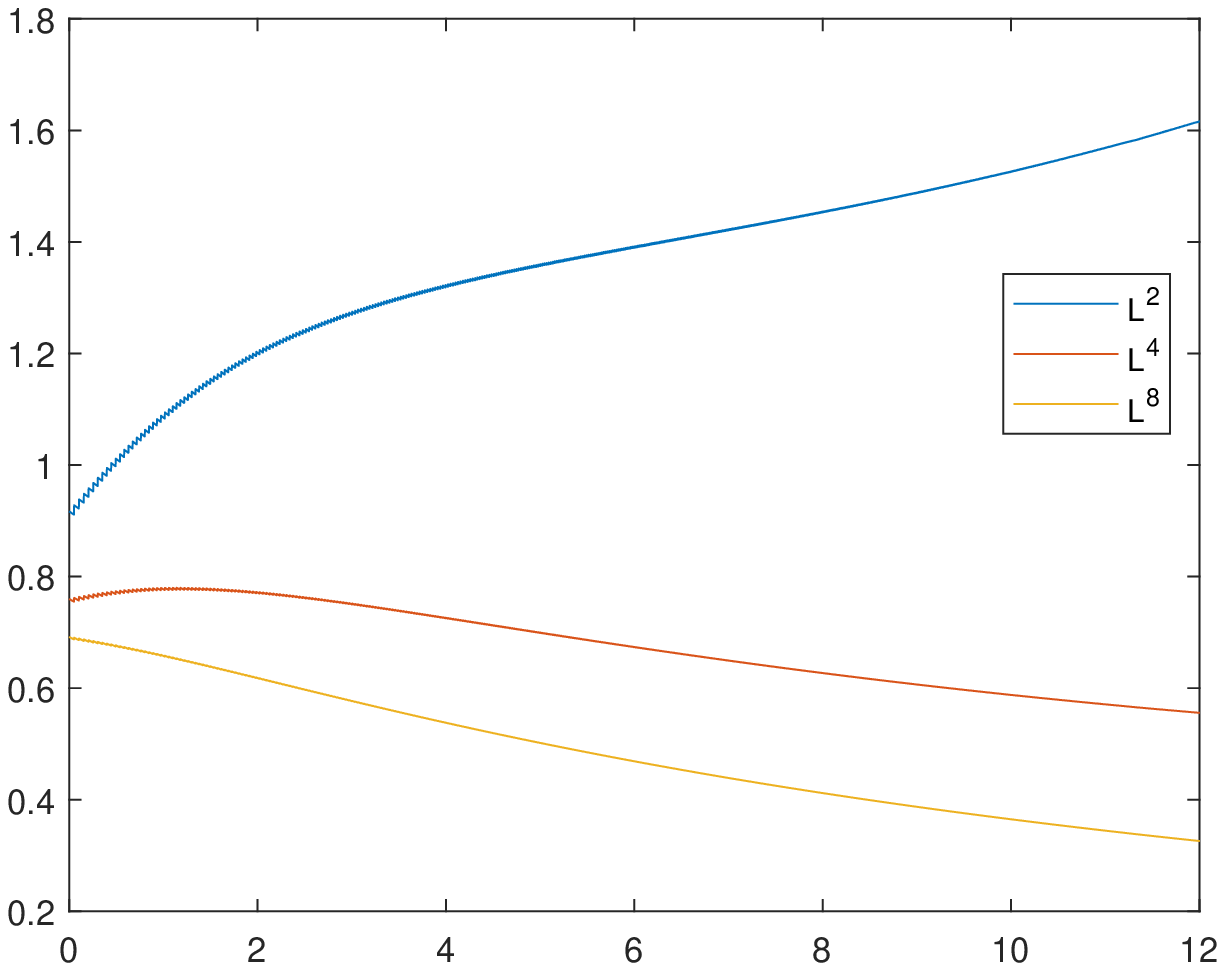}
\caption{$D=0.1,K_1=0.1,K_2=0.01$}
		\label{fig:l2n2}
\end{subfigure}
\caption{Evolution of $||\mu(\cdot,t)-\mu^*||_{L^{2n}({\Omega(t)})}$ with respect to time under different parameter regimes. Parameters: $\gamma=80$, $G(c)=gc,\psi(c)=c$ with $g=1,a=0.5$ and $c_B=1$. Here the threshold value of nutrient concentration $c_0=0.5$.} 
\label{fig:s3l2n2}
	\end{figure}

In Figure \ref{fig:s3l2n2} we numerically simulate the compressible model to illustrate the competition between exponential growth of $|\Omega(t)|$ and exponential decay of $|\mu-\mu^*|$. We plot the evolution of $||\mu(\cdot,t)-\mu^*||_{L^{2n}({\Omega(t)})}$ with different $n$ under two parameters regimes. In the first regime (Figure \ref{fig:l2n1}) $K_2$ is large and the $L^2,L^4,L^8$ norm decay fast. In the second regime (Figure \ref{fig:l2n2}) $K_2$ is small. $L^2$ norm tends to diverge to infinity but $L^4,L^8$ norm are decreasing. This also indicates that the condition \eqref{eq:condition-thm2} in Theorem \ref{thm:L2n} may not be sharp since it is violated for all the three cases $L^2,L^4,L^8$ in the second regime.

The ``well-mixed'' limit shows that under Assumption \ref{as:non-negative}, the density fraction $\mu(x,t)$ will converge to a spatial-homogeneous constant state exponentially fast. This allows us to further simplify the model by taking $\mu$ to be a constant for all $x\in\overline{\Omega(t)}$, i.e. $\mu\equiv\mu^*$, where $\mu^*$ is the root in $(0,1)$ of $f$ in \eqref{eq:rec-ode}. For this simplified model with constant $\mu$, it is convenient to draw analytical solutions as in \cite{liu2019analysis} and to further study the behavior of the free boundary model, as we do in the next section.


\section{Analytical solution and influence of autophagy on expansion speed}\label{sc:auto}

In this section, we first derive an analytical solution for the free boundary model with constant $\mu$, thanks to the well-mixed limit in section 3 Theorem \ref{thm:uni-con}. From the analytical solution we can analyze the influence of autophagy. In particular, it shows that when the nutrient supply rate by autophagic cells $a$ is relatively larger than the extra death rate $D$ from self-killing, autophagy will result in exponential growth rate of the tumor, in contrast to the linear growth rate in a related tumor-nutrient model in \cite{liu2019analysis}. The essential mechanism is that in the presence of the autophagy, sufficient nutrients are available to the tumor, whatever big it is. Also, we do numerical simulations on the compressible cell density model for demonstration of analytical results and for exploration of more scenarios.

\subsection{A one-dimensional analytic solution with constant $\mu$}

The ``well-mixed'' limit in Section \ref{sc:math-proof} shows that under Assumption \ref{as:constants}, the density fraction $\mu$ has a constant steady state: $\mu\equiv \mu^*$, for all $\ x\in\Omega(t)$. Here $\mu^*$ as in \eqref{eq:roots} is the root of $f$ \eqref{eq:rec-ode}. This allows use to simplify the model by set $\mu$ to be the constant $\mu^*$, i.e., $\mu\equiv\mu^*$, $\forall x\in\Omega(t),t\geq0$. Then the model \eqref{eq:free-init}$\sim$\eqref{eq:model-v} reduces to the evolution of $(\Omega,p,c)$ \eqref{eq:model-p}$\sim$\eqref{eq:model-v}:
\begin{equation}\label{model:reducedPDE}
\begin{cases}
-\Delta p=\mu G(c)+(1-\mu)(G(c)-D),\quad x\in\Omega(t).\\
p=0,\quad x\in \partial\Omega(t).\\
-\Delta c+\psi(c)=(1-\mu)a  ,\quad x\in \Omega(t),\\
c=c_B>0  ,\quad  x\in\partial\Omega(t),\\
V_n=-\nabla p \cdot \textbf{n},\quad x\in\p\Omega(t).
\end{cases}
\end{equation}

To derive an analytical solution, we shall specify the parameters. For the net growth rate of normal cells $G$ we assume the proliferation rate is proportional to the concentration of nutrients $c$ with a factor $g$ and we omit the death rate of normal cells for simplicity. Thus $G(c)=gc$. For the consumption rate $\psi(c)$, we assume it is proportional to $c$. Precisely we choose $\psi(c)=c$. And since nutrient is relatively sufficient at the boundary of tumor, we assume the threshold value for autophagy cells $c_0$ in \eqref{ap-3psi} is less than the boundary value for the concentration of nutrients $c_B$. Note when $\psi(c)=c$, $c_0=a$. To summarize, we assume
\begin{equation}\label{para:sc41}
G(c)=gc,\quad\psi(c)=c,\quad a<c_B.
\end{equation}

Furthermore, for simplicity we consider one-dimensional symmetric case, when the tumor region $\Omega(t)=(-R(t),R(t))$. With these choices of parameters, \eqref{model:reducedPDE} becomes
\begin{equation}\label{model:analy}
\begin{cases}
-\frac{d^2}{dx^2}p=gc-(1-\mu)D,\quad  x\in (-R(t),R(t)),\\
p=0,\quad x=\pm R(t).\\
-\frac{d^2}{dx^2}c+ c=(1-\mu)a,\quad  x\in (-R(t),R(t)),\\
c=c_B>0,\quad x=\pm R(t).\\
\frac{dR}{dt}=-\frac{dp}{dx}|_{x=R(t)}.
\end{cases}
\end{equation}System \eqref{model:analy} is essentially a single-species model. And the following derivation of an analytical solution is similar to that in \cite{liu2019analysis}.

First we solve $c$ from the equation of nutrients. With the boundary condition $c(R(t))=c_B$ and the symmetric condition $\frac{d}{dx}c|_{x=0}=0$, we obtain the expression of $c$:
\begin{align}
c(x,t)&=(1-\mu)a+\frac{c_B-(1-\mu)a}{\cosh(R(t))}\cosh(x)\label{cnew}\\
&=\frac{c_B}{\cosh(R(t))}\cosh(x)+(1-\mu)a(1-\frac{\cosh(x)}{\cosh(R(t))}). \notag
\end{align}
\begin{figure}[H]
	\centering
	\begin{subfigure}{.49\textwidth}
		\centering
		\includegraphics[width=0.9\textwidth]{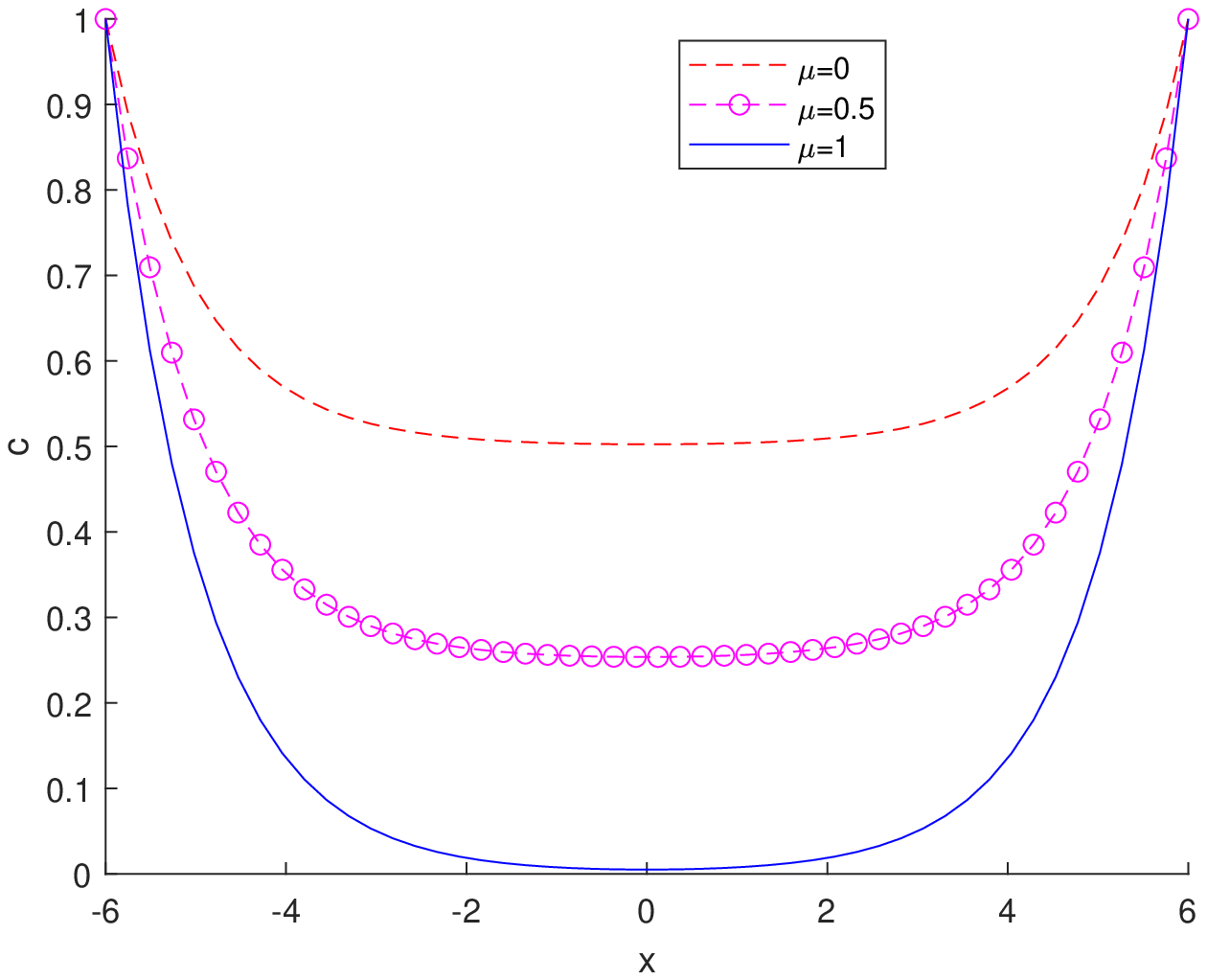}
		\caption{Concentration of nutrients $c$}
		\label{g-c}
	\end{subfigure}
	\begin{subfigure}{.49\textwidth}
		\centering
		\includegraphics[width=0.9\textwidth]{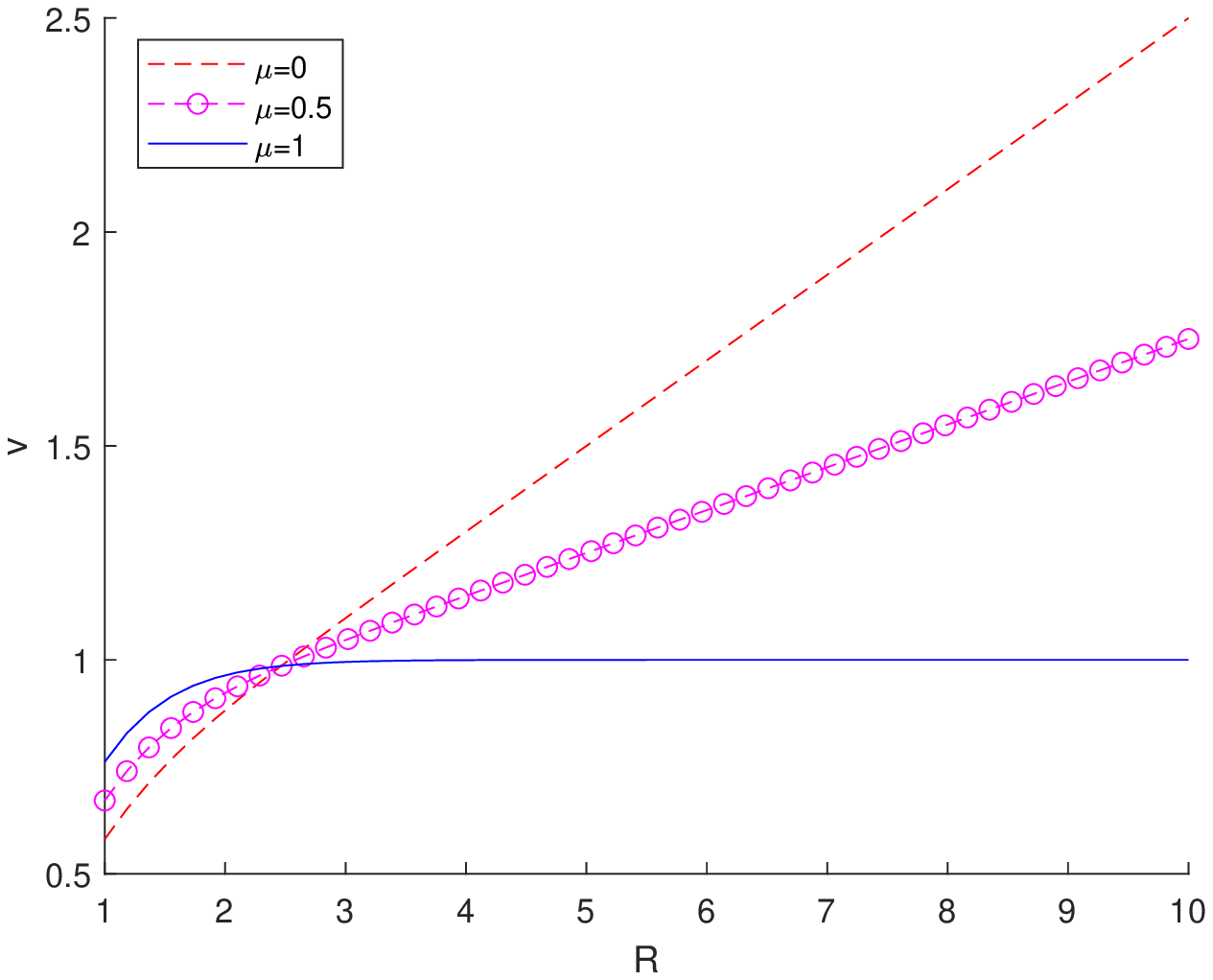}
		\caption{Moving speed of the boundary $v:=\frac{dR(t)}{dt}$ w.r.t $R$}
		\label{g-r}
	\end{subfigure}
	\caption{Graph of nutrient concentration $c$ and moving speed of the boundary for different $\mu$. When $\mu=1$ there is no autophagic cells, when $\mu=0$ all cells are autophagic cells. In the presence of autophagy, the nutrient is more sufficient and the growth rate tends to be linear in $R$, which leads to an exponential growth of $R$. Parameters: $c_B=1,a=0.5,g=1,D=0.3$.
	}
	\label{fig:analy-cv}
\end{figure}
Plugging \eqref{cnew} into the equation for $p$, one gets
\begin{align*}
-\frac{d^2}{dx^2}p&=gc-(1-\mu)D=\frac{g(c_B-(1-\mu)a)}{\cosh(R(t))}\cosh(x)+(1-\mu)ga-(1-\mu)D.
\end{align*}
Together with boundary condition $p(\pm R(t))=0$, the solution of $p$ is
\begin{equation}\label{eq:ana-p}
p(x,t)=-\frac{g(c_B-(1-\mu)a)}{\cosh(R(t))}\cosh(x)-\frac{1}{2}(1-\mu)(ga-D)x^2+C,
\end{equation}
where the coefficient of linear term $x$ is $0$ by symmetry and $C$ could be determined from $p(R(t))=0$. Actually 
$C=g(c_B-(1-\mu)a)+\frac{1}{2}(1-\mu)(ga-D)R(t)^2$. 

In Figure \ref{g-c} we plot the graph of nutrients $c$ with different $\mu$. When all cells are normal cells ($\mu=1$), in the middle of tumor the concentration of nutrients $c$ decays to zero, which indicates that the proliferation is limited. While in the presence of autophagic cells ($\mu=0.5,1$), the concentration of nutrients can maintain a basal level to support proliferation in the core of tumor.

From \eqref{eq:ana-p} the moving speed of the boundary reads:
\begin{align}\label{v1}
\frac{dR(t)}{dt}=-\frac{dp}{dx}|_{x=R(t)}&=g(c_B-(1-\mu)a)\tanh(R(t))+(1-\mu)(ga-D)R(t).
\end{align} 

Whereas, when there is no autophagic cells (i.e., $\mu=1$), the moving speed is, as that in \cite{liu2019analysis}:
\begin{equation}\label{v2}
\frac{dR(t)}{dt}=gc_B \tanh(R(t)).
\end{equation}

\begin{figure}[H]
\centering
	\includegraphics[width=0.6\textwidth]{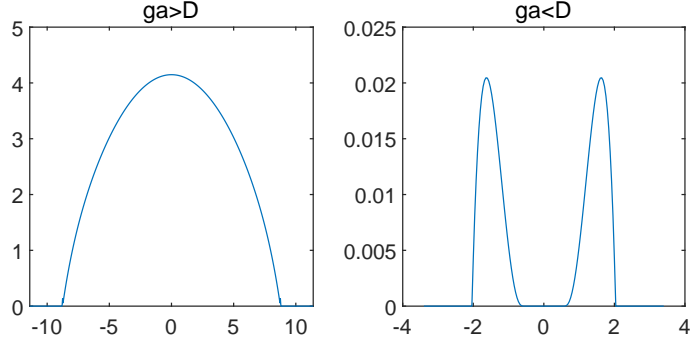}
	\caption{\hl{Graph of pressure $p$ for cell density model for $\gamma=80$. Left: $ga>D$. Right $ga<D$, there is a necrotic core in the middle.}}
	\label{fig:p}
\end{figure}

\hl{In \eqref{v2} as $R\rightarrow\infty$, the moving speed $\frac{dR}{dt}\rightarrow gc_B$, a constant.} An essential difference brought by autophagy is the linear term on the right hand side of \eqref{v1}, compared to \eqref{v2}. Recall our choice of parameters \eqref{para:sc41}, a dichotomy arises, which depends on the growth factor $g$, the nutrient supply rate by autophagy $a$ and the extra death rate due to autophagy $D$:
 \begin{enumerate}
 	\item  When $ga>D$, the linear term will result in exponential growth, while when there is no autophagy the growth rate in \eqref{v2} is nearly a constant, as plotted in Figure \ref{g-r}, and results in linear growth. Biologically this indicates that when the nutrients supplied by autophagic cells are sufficient to offset its extra death, autophagy will significantly accelerate tumor growth. \hl{Together with the nutrient concentration shown in Figure \ref{g-c}, this is in accordance with the biological observation that autophagy helps cells in the center of tumor, where the nutrient is lack, and therefore promotes tumor growth \cite{crwhite2012deconvoluting,kang2008or,degenhardt2006autophagy}.}
 	\item While, on the contrary when $ga<D$, the solution may violate the non-negative Assumption \ref{as:non-negative} if $R(t)$ is large. This indicates when the extra death of autophagy $D$ is dominant, autophagy may accelerate cells death and contributes formation of necrotic core. \hl{In this case, the former description \eqref{model:reducedPDE} or \eqref{model:analy} of free boundary problem is no longer correct for the incompressible limit $\gamma\rightarrow\infty$ of the cell density model \eqref{model:PDE-systems}. Although the necrotic core is not the focus of this paper, a simulation of cell density model \eqref{model:PDE-systems} with large $\gamma$ is given in Figure \ref{fig:p} for illustration. Actually for a correct free boundary description of the limit solution, the Poisson equation of pressure $p$ in \eqref{model:reducedPDE} or \eqref{model:analy} should be replaced by the corresponding obstacle problem in $\Omega(t)$ with the constraint $p\geq0$. This dichotomy on the sign of $ga-D$ partially reflects the dual role of autophagy -- both pro-survival and pro-death \cite{lin2015autophagy}.}
 \end{enumerate}

We could understand the exponential growth in the presence of autophagy in a more direct way. From the conversation of density, we have
\begin{equation}\label{eq:speed-conver}
\frac{dR(t)}{dt}=\frac{1}{2}\int_{-R(t)}^{R(t)}(G(c)-(1-\mu)D) dx.
\end{equation} 
By the maximum principle we have a lower bound for $c$, precisely $c\geq\min\{c_B,(1-\mu)a\}=(1-\mu)a$. Thus we have $G(c)\geq(1-\mu)ga$. Combine this with \eqref{eq:speed-conver}, we obtain
\begin{equation*}
\frac{dR(t)}{dt}=\frac{1}{2}\int_{-R(t)}^{R(t)}(G(c)-(1-\mu)D) dx\geq \frac{1}{2}\int_{-R(t)}^{R(t)}(1-\mu)(ga-D)dx=(1-\mu)(ga-D)R(t).
\end{equation*}Thus we get the exponential growth of $R$:
\begin{equation}\label{eq:exp-growth}
R(t)\geq \exp{\{(1-\mu)(ga-D)t\}}R(0).
\end{equation}This reasoning can easily extend to higher dimensional cases with similar assumptions.
\subsection{Numerical simulations}
In this subsection, we numerically simulate the compressible PDE model \eqref{model:PDE-systems} with large $\gamma$ to demonstrate the incompressible limit and analysis on the moving speed of the boundary.

For the numerical scheme, we adapt the prediction-correction framework in \cite{liu2018accurate} which we present in details in Appendix \ref{app:num}. Our computations are in one-dimension. We choose the spatial step $\Delta x=0.04$ and the temporal step $\Delta t=0.002$.

We recall in the compressible model, densities of normal cells $n_1$ and autophagic cells $n_2$ satisfy the following system:
\begin{equation*}
\begin{cases}
\frac{\partial n_1}{\partial t}-\Div(n_1 \nabla p)=G(c)n_1-K_1n_1+K_2n_2,\\
\frac{\partial n_2}{\partial t}-\Div(n_2 \nabla p)=(G(c)-D)n_2+K_1n_1-K_2n_2,\\
p=\frac{\gamma}{\gamma-1}(n_1+n_2)^{\gamma-1},\\
\end{cases}
\end{equation*}We choose compactly supported total density $n=n_1+n_2$ for initial data. Thus $\Omega(t):=\{x:n(x,t)>0\}$ is a bounded domain for all $t>0$, which allows us to choose the quasi-static case of \eqref{eq:c-compact} for the evolution of nutrients:
\begin{equation*}
\begin{cases}
-\Delta c+\psi(c)(n_1+n_2)=an_2,\quad x\in\Omega(t),\\
c=c_B,\quad x\notin\p\Omega(t).
\end{cases}
\end{equation*}And parameters are chosen of forms in  \eqref{para:sc41}.

\begin{figure}[H]
	\centering
\begin{subfigure}{.49\textwidth}
	\includegraphics[width=0.9\textwidth]{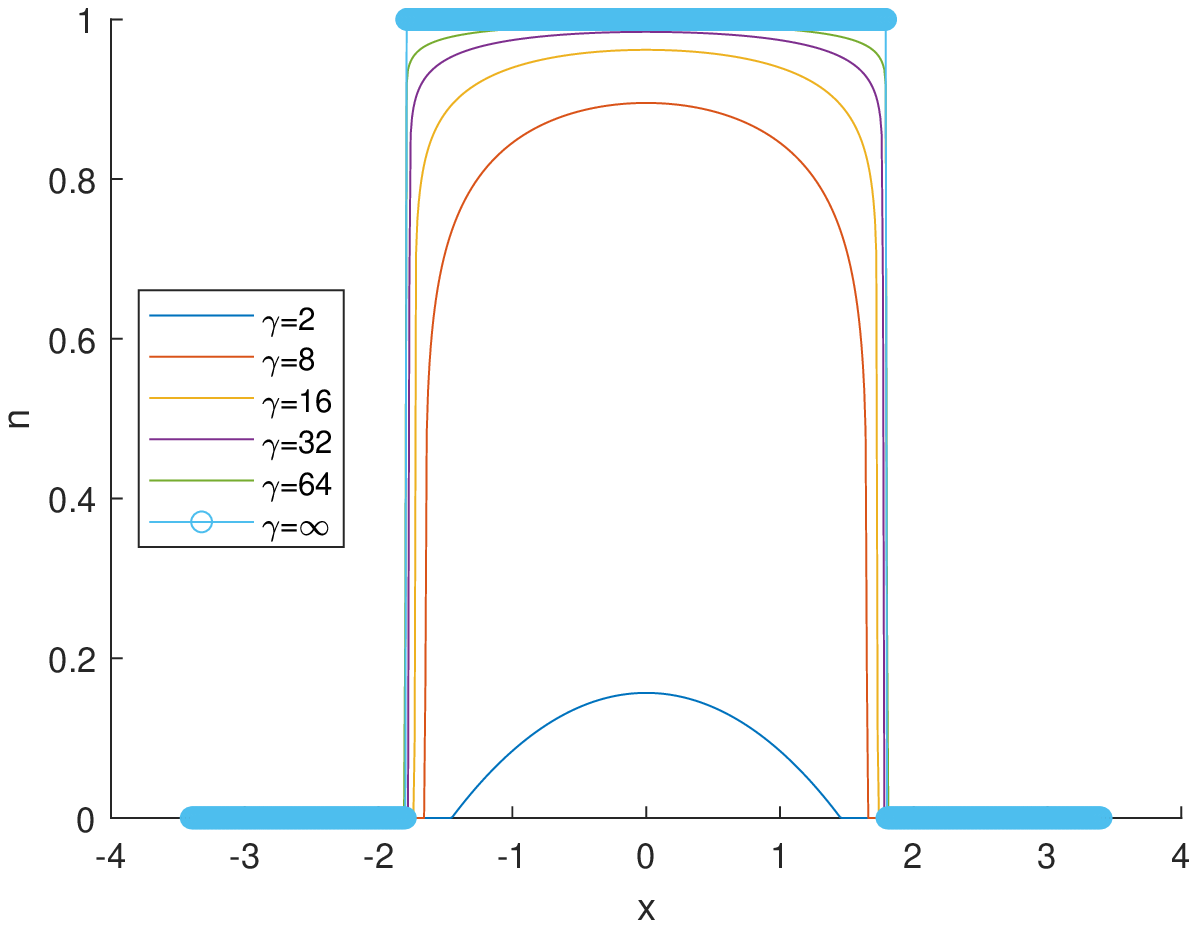}
	\caption{Total cell density $n$}
\end{subfigure}
\begin{subfigure}{.49\textwidth}
	\includegraphics[width=0.9\textwidth]{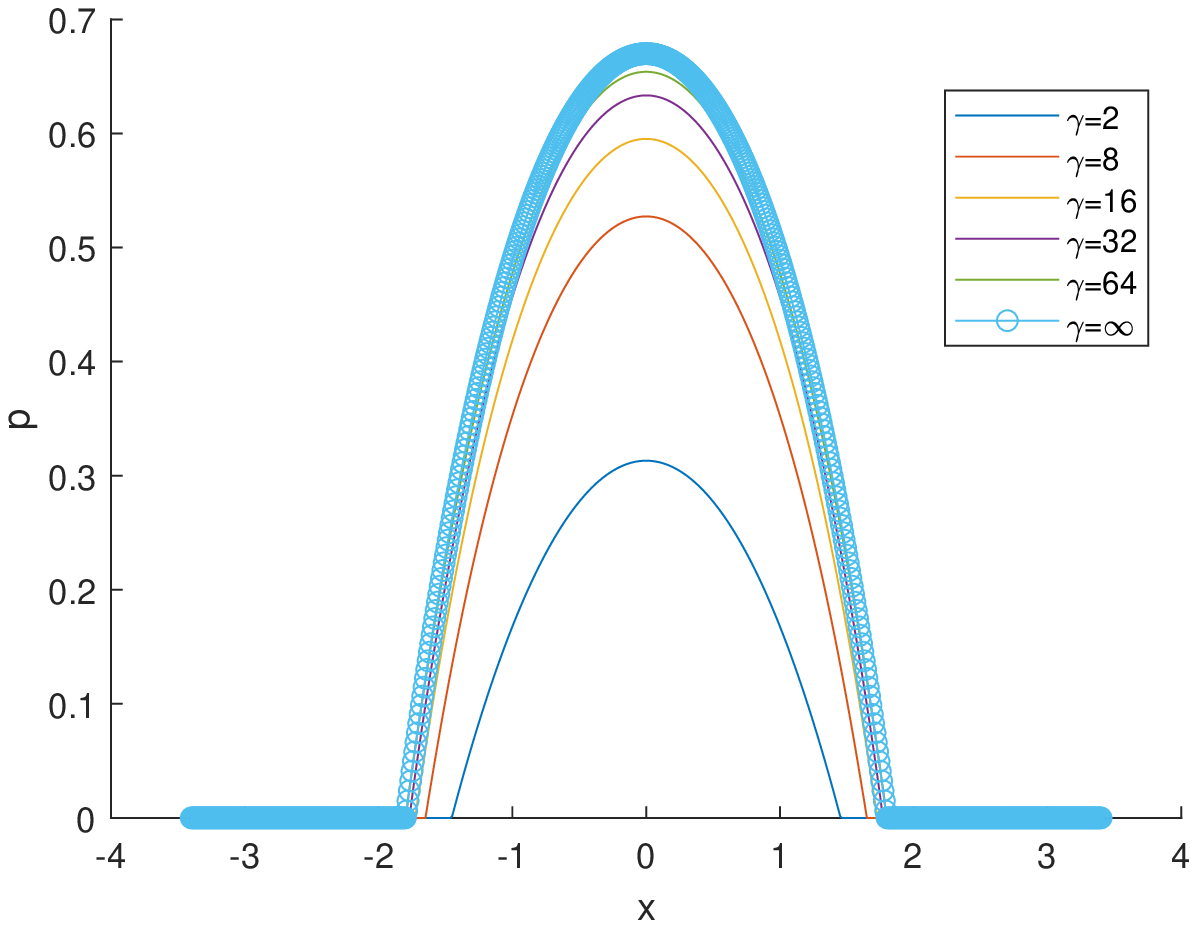}
	\caption{Pressure $p$}
\end{subfigure}
\caption{Plots of the total density $n$ and the pressure $p$ for different $\gamma$ at $t=1$. Left: $n$. Right $p$. $\gamma=\infty$ stands for the analytical solution of the free boundary problem. As $\gamma$ increases, the solution of the compressible model approximates the solution of the free boundary model. Parameters: $g=1,a=0.5,D=0.3,c_B=1,K_1=K_2=1$}
\label{fig:s4limit}
\end{figure}

\begin{figure}[H]
	\centering
	\begin{subfigure}{.49\textwidth}
		\centering
		\includegraphics[width=0.8\textwidth]{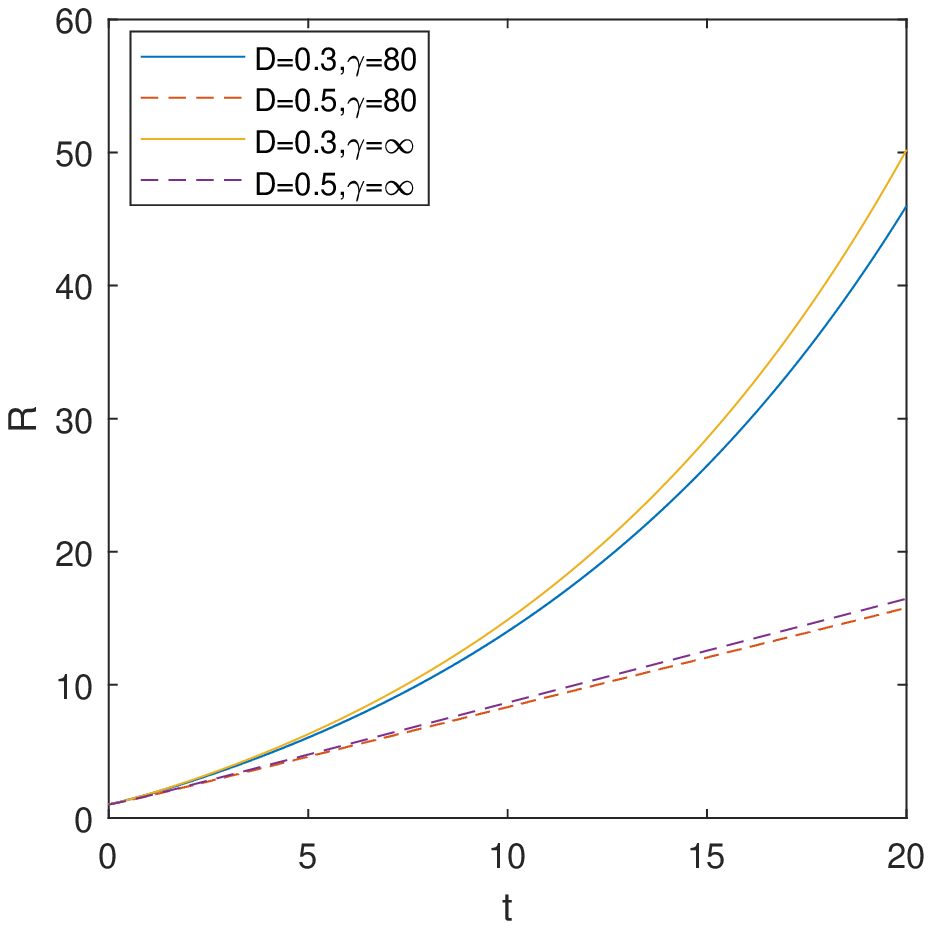}
		\caption{Tumor radius $R(t)$}
		\label{fig:s4f21}
	\end{subfigure}
	\begin{subfigure}{.49\textwidth}
		\centering
		\includegraphics[width=0.8\textwidth]{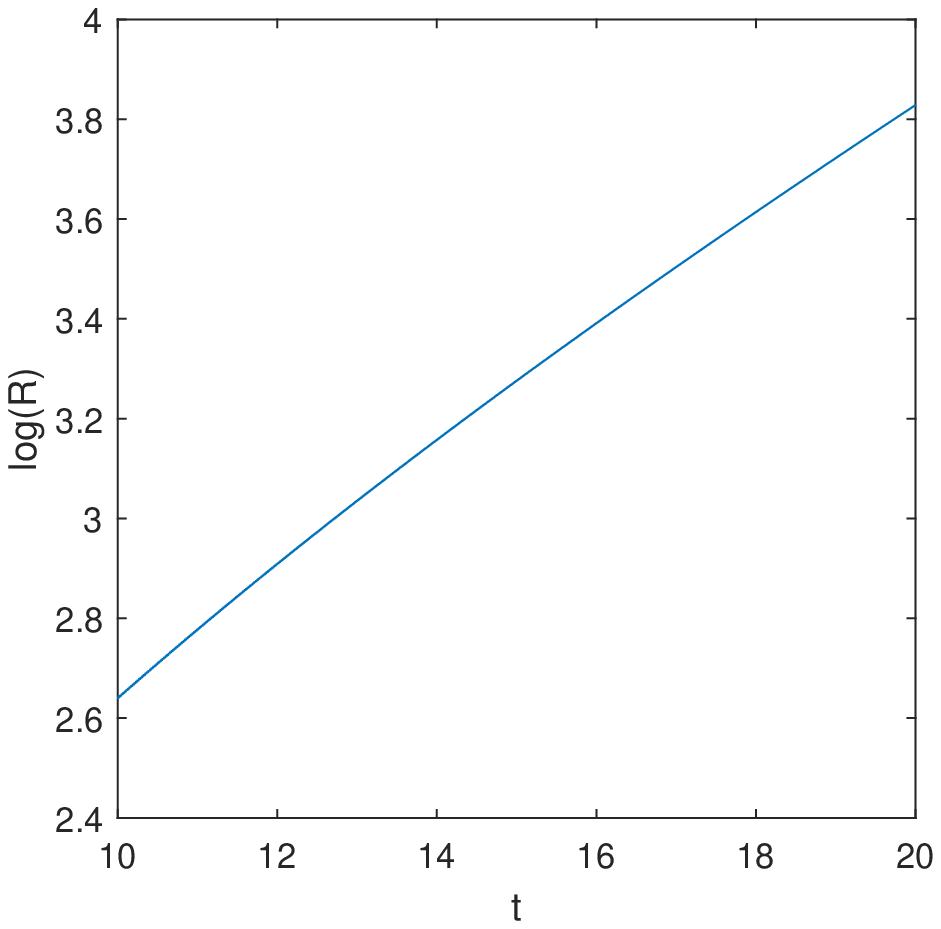}
		\caption{$\log(R(t))$  for $D=0.3$.}
				\label{fig:s4f22}
	\end{subfigure}
	\caption{Growth of tumor radius w.r.t time for different $D$. Other parameters $g=1,a=0.5$. Left: radius w.r.t. time. Solid line: $D=0.3$ and thus {{$ga>D$}}. Dashed line: $D=0.5$ and thus {$ga=D$}. Right: $\log R$ w.r.t. time for $D=0.3$ from $t=10$ to $t=20$.}
	\label{fig:s4f2}
 \end{figure}

First, to demonstrate the incompressible limit, we compare numerical solutions of total density $n$ and pressure $p$ for different $\gamma$ with analytical solutions of the free boundary model (denoted by $\gamma=\infty$ formally). Initially we take $\Omega(0)=(-R,R)$ with $R=1$ and specify the pressure $p$ according to the analytical solution of free boundary model \eqref{eq:ana-p}. Then the initial total density $n$ is recovered from \eqref{eq:intro-p}. We choose initial density fraction $\mu$ to be at the well-mixed constant $\mu^*$, thus we could get initial density of two kinds of cells $n_1=\mu^*n$ and $n_2=(1-\mu^*)n$. For different $\gamma$ we numerically evolves this system till time $t=1$. We plot these numerical solutions as well as the analytical solution. For total density $n$ the analytical solution is $\mathbb{I}_{(-R(t),R(t))}$ while for pressure $p$ that is given in \eqref{eq:ana-p}. Here $R(t)$ is obtained by numerically solve the ODE \eqref{v1}. Results in Figure \ref{fig:s4limit} show that as $\gamma$ increases the numerical solutions get closer to the analytical solution of the free boundary model.
Next, we investigate the boundary propagation speed, which are to be compared with the analytical results. For compressible models with large $\gamma$, we plot the evolution of radius with respect to time in two cases: $ga=D$ and $ga>D$ in Figure \ref{fig:s4f21}. Results are close to the evolution of $R(t)$ in the free boundary model, i.e., the ODE \eqref{v1}. We also plot $\log R$ with respect to time in the case $ga>D$ to verify the exponential growth in Figure \ref{fig:s4f22}. Numerical solutions are consistent with the analytical results: when $ga=D$ the growth is nearly linear and when $ga>D$ the radius tends to grow exponentially.

\begin{figure}[H]
	\centering
	\begin{subfigure}{\textwidth}
		\centering
		\includegraphics[width=0.8\textwidth]{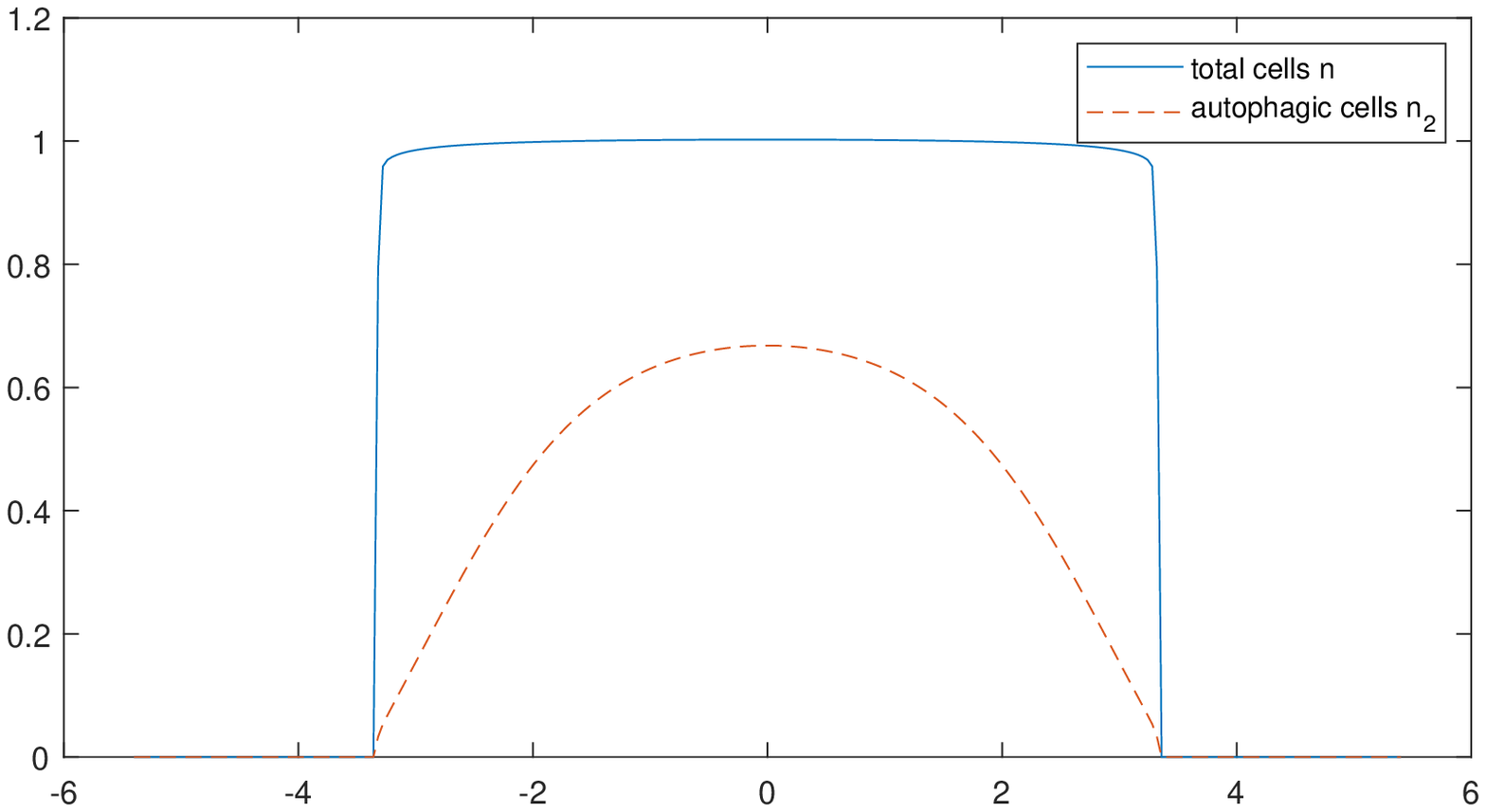}
		\caption{\hl{Graph of total cell density $n$, autophagic cell density $n_2$. Dashed line: $n_2$. Solid line: $n=n_1+n_2$.}}
		\label{fig:s4fin1}
	\end{subfigure}
	\begin{subfigure}{\textwidth}
		\centering
	\includegraphics[width=0.8\linewidth]{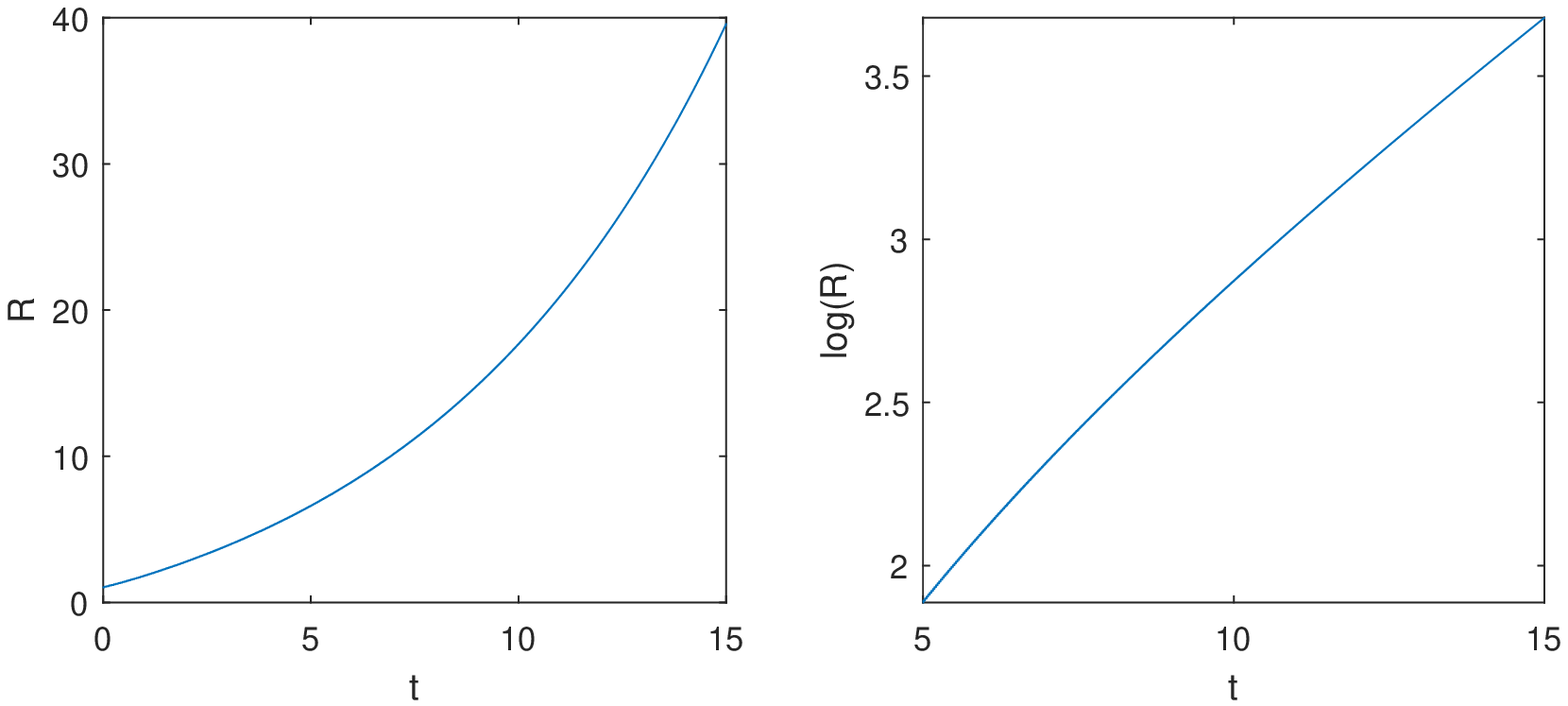}
	\caption{\hl{Growth of the tumor radius $R$}}
			\label{fig:s4finR}
	\end{subfigure}
\caption{Plots of the total cell density $n$ and the autophagic cell density $n_2$, and the evolution of tumor radius in the case $K_1,K_2$ are not constants. Parameters $a=0.5,D=0.3,\gamma=80,K_1(c)=(\frac{1-c}{c+0.1})_+,K_2(c)=\frac{2c}{c+1}$.}

\end{figure}

\hl{
Then, we consider the case $K_1,K_2$ are not constants, i.e., Assumption \ref{as:constants} does not hold. In this case, $\mu=\frac{n_1}{n_1+n_2}$ needs not to be a constant and can vary spatially. Considering the case $ga>D$, we plot the graph of total cell density $n$ and autophagic cell density $n_2$ in Figure \ref{fig:s4fin1}. The ratio of autophagic cells is high in the center of tumor, where the nutrient is lack, and is small at the tumor boundary, where the nutrient is sufficient. This captures the experiment observation that autophagy prominently locates in tumor cells under deprivation of nutrients \cite{degenhardt2006autophagy,crwhite2012deconvoluting}. Also we plot the evolution of radius, which shows the tumor still tends to grow exponentially in this case.
}

\section*{Acknowledgements}
 X. Dou is partially supported by the elite undergraduate training program of School of Mathematical Sciences in Peking University. J. Liu is partially supported by National Science Foundation (NSF) under award DMS-2106988. Z. Zhou is supported by the National Key R\&D Program of China, Project Number 2020YFA0712000 and NSFC grant No. 11801016, No. 12031013. The authors thank Jinzhi Lei and Beno{\^\i}t Perthame for helpful discussions.

	
	
	
    

\begin{appendices}
	
\section{The proof of Proposition 2 }\label{app:proof}

\begin{proof}[Proof of Proposition \ref{prop:refor}]
	Since $\nabla p$ is Lipschitz in space thanks to $p$ is $C^2$ in space by Assumption \ref{as:strong-simple-con}, we could applies Picard's theorem to obtain the uniqueness and local existence. 

	Note that when $y\in\partial\Omega(0)$ the ODE is exactly the evolution of boundary in Assumption \ref{as:strong-simple-con}. Thus it is well-defined and by our definition $\p\Omega(t)=X_{t,0}(\p\Omega(0))$.

	Then we could obtain $x(t,y)\in\Omega(t)$ (if the flow exists for time $t$) for $y\in\Omega(0)$. Combine this with continuity of $\nabla p$, we could extend the local solution to get the existence in $[0,T)$.

	It remains to show  $X_{t,0}(\Omega(0))= \Omega(t)$,. We have obtained  $X_{t,0}(\Omega(0))\subseteq\Omega(t)$, it remains to prove that for all $x\in\Omega(t)$, there exists $y\in\Omega(0)$ such that $x=X_{t,0}(y)$. This can be obtained by solving flow map (\ref{eq:flowmap}) backward and apply the same argument in the first part of proof.
\end{proof}
%

	\section{Numerical scheme for the compressible model}\label{app:num}
	We adapt the numerical scheme in \cite{liu2018accurate} for the compressible cell density model \eqref{model:PDE-systems}. In this section we present the details of the scheme. For properties and analysis on the scheme we refer to \cite{liu2018accurate}.
		
	The main difficulty is the nonlinearity in \eqref{eq:intro-p} with high $\gamma$. In the projection and correlation framework of \cite{liu2019analysis}, they introduce the velocity field $\textbf{u}$:
	\begin{equation*}
	\textbf{u}:=-\nabla p.
	\end{equation*} and solve the equation for $\textbf{u}$,
	\begin{equation}\label{eq:u}
	\frac{\p \textbf{u}}{\p t}=\gamma\nabla(n^{\gamma-2}(\nabla\cdot(n\textbf{u})-n_1G_1(c)-n_2G_2(c))).
	\end{equation}Here we use the shorthand notation $G_1(c)=G(c)$, $G_2(c)=G(c)-D$.
	
	We describe the update from $(u^j,n_1^j,n_2^j,n^j,c^j)$ to $(u^{j+1},n_1^{j+1},n_2^{j+1},n^{j+1},c^{j+1})$ in one-dimension.
	
	Consider the domain $\Omega=[a,b]$. $\Delta x=\frac{b-a}{N}$ be the spatial mesh size, and the grid points are
	\begin{equation*}
	x_i=a+i\Delta x,\ x_{i+\frac{1}{2}}=a+(i+\frac{1}{2})\Delta x,
	\end{equation*}in regular grid and staggered grid, respectively. As in \cite{liu2018accurate} the regular grid is for total density $n$ and the staggered grid is for the velocity field $u$.

	First we solve the equation for velocity field $u$ \eqref{eq:u} to get $u^*$:
	
	\begin{equation}
	\begin{aligned}
		\frac{u_{i+1 / 2}^{j *}-u_{i+1 / 2}^{j}}{\Delta t}=& \frac{\gamma}{\Delta x}\left\{\left(\left(n_{i+1}^{j}\right)^{\gamma-2}\left(\frac{n_{i+3 / 2}^{j} u_{i+3 / 2}^{j *}-n_{i+1 / 2}^{j} u_{i+1 / 2}^{j *}}{\Delta x}-n_{1,i+1}^{j} G_{1,i+1}^{j}-n_{2,i+1}^{j} G_{2,i+1}^{j}\right)\right)\right.\\
		&\left.-\left(\left(n_{i}^{j}\right)^{\gamma-2}\left(\frac{n_{i+1 / 2}^{j} u_{i+1 / 2}^{j *}-n_{i-1 / 2}^{j} u_{i-1 / 2}^{j *}}{\Delta x}-n_{1,i}^{j} G_{1,i}^{j}-n_{2,i}^{j} G_{2,i}^{j}\right)\right)\right\},
	\end{aligned}
	\end{equation}
	Here $n_{1,i},n_{2,i}$ are the spatial discretization for $n_1$, $n_2$. And $G_{1,i}^j=G_1(c_{i}^j)$, $G_{2,i}^j=G_2(c_{i}^j)$. The value of $n$ on the staggered grid is approximated by
	\begin{equation*}
	n_{i+\frac{1}{2}}^j=\frac{n^j_{i}+n^{j}_{i+1}}{2}.
	\end{equation*}

	With $u^*$ obtained, the two equations for $n_1,n_2$ are discretized by a central scheme, and we treat the reaction term semi-implicitly.
	\begin{equation}
	\begin{aligned}
	\frac{n_{1,i}^{j+1}-n_{1,i}^{j}}{\Delta t}+\frac{F_{1,i+1 / 2}^{j}-F_{1,i-1 / 2}^{j}}{\Delta x}=G_{1,i}^{j}n_{1,i}^{j+1}-K_{1,i}^jn_{1,i}^{j+1}+K_{2,i}^jn_{2,i}^{j+1}\\
	\frac{n_{2,i}^{j+1}-n_{2,i}^{j}}{\Delta t}+\frac{F_{2,i+1 / 2}^{j}-F_{2,i-1 / 2}^{j}}{\Delta x}=G_{2,i}^{j}n_{2,i}^{j+1}+K_{1,i}^jn_{1,i}^{j+1}-K_{2,i}^jn_{2,i}^{j+1}.
	\end{aligned}\end{equation}Here $K_{1,i}^j=K_1(c_{i}^j)$, $K_{2,i}^j=K_2(c_{i}^j)$.

	Now we describe how to compute the flux $F_1,F_2$. For clarity we omit the dependence on the specific kind of cells and, with a bit abuse of notation, use $n_{i}^j$ instead of $n_{1,i}^j,n_{2,i}^j$. When computing $F_1$, $n_i^j$ in the following should be substituted by $n_{1,i}^j$. Similarly when computing $F_2$, $n_i^j$ should be substituted by $n_{2,i}^j$.
	\begin{equation*}
	F_{i \pm 1 / 2}^{j}=\frac{1}{2}\left[n^{L j} u^{j *}+n^{R j} u^{j *}-\left|u^{j *}\right|\left(n^{R j}-n^{L j}\right)\right]_{i \pm 1 / 2}
	\end{equation*}Here the edge values $n_{i\pm 1/2}^{L/Rj}$ are defined as follows:
	\begin{equation*}
	n_{i+1/2}^{Lj}=n_{i}^j+\frac{\Delta x}{2}(\p_xn)_i^j,\quad n_{i+1/2}^{Rj}=n_{i+1}^j-\frac{\Delta x}{2}(\p_xn)_{i+1}^j,
	\end{equation*}where $(\p_xn)_{i}^j$ is given by
	\begin{equation*}
	(\p_xn)_{i}^j=\begin{cases}
	\min\{\frac{n_{i+1}^j-n_{i}^j}{\Delta x},\frac{n_{i+1}^j-n_{i-1}^j}{2\Delta x},\frac{n_{i}^j-n_{i-1}^j}{\Delta x}\},\quad \text{if all are positive},\\	\max\{\frac{n_{i+1}^j-n_{i}^j}{\Delta x},\frac{n_{i+1}^j-n_{i-1}^j}{2\Delta x},\frac{n_{i}^j-n_{i-1}^j}{\Delta x}\},\quad \text{if all are negative},\\	0,\quad \text{otherwise}.\\
	\end{cases}
	\end{equation*}
	After obtaining $n_{1,i}^{j+1}$ and $n_{2,i}^{j+1}$, we compute the new total density $n$ and new velocity field $u$:
\begin{equation}
n_{i}^{j+1}=n_{1,i}^{j+1}+n_{2,i}^{j+1},\quad u_{i+1/2}^{j+1}=-\frac{\gamma}{\gamma-1}\frac{(n_{i+1}^{j+1})^{\gamma}-(n_{i}^{j+1})^{\gamma}}{\Delta x}.
\end{equation}
	Finally,  for the concentration of nutrients $c$, since we simulate the linear case $\psi(c)=c$, we use standard finite difference method for the diffusion equation with Neumann boundary condition:
	\begin{equation}
	\begin{aligned}
		\frac{c_i^{j+1}-c_{i}^{j}}{\Delta t}-\frac{c_{i+1}^{j+1}-2c_{i}^{j+1}+c_{i-1}^{j+1}}{(\Delta x)^2}+c_{i}^{j+1}n_{i}^{j+1}=an_{2,i}^{j+1},\\\frac{c_1^{j+1}-c_0^{j+1}}{\Delta x}=\frac{c_{N-1}^{j+1}-c_{N}^{j+1}}{\Delta x}=\lambda^{j+1}=\lambda((j+1)\Delta t).
	\end{aligned}
	\end{equation}

	For the boundary condition, by choosing the computation domain larger than the support of $n$ we reduce the problem to a Dirichlet problem with zero boundary condition. In simulation we track the support of $n$ and dynamically enlarge the computation domain such that it is much larger than the support of $n$. In the support of $n$ we solve the equation of $c$ \eqref{eq:c-compact} by standard finite difference method.
		

	
\end{appendices}
\bibliography{intro_arxiv_2021}
\bibliographystyle{plain}
\end{document}